\documentclass[a4paper,12pt]{article}
\usepackage{amsmath}
\usepackage{amsfonts}
\usepackage{amssymb}
\usepackage{latexsym}
\usepackage{epsfig}
\usepackage{graphicx}
\usepackage{oldgerm}
\usepackage{theorem}

\def\C{\mathbb{C}}
\def\E{{\mathbb E}}
\def\mbK{\mathbb{K}}
\def\N{\mathbb{N}}
\def\P{{\mathbb P}}

\def\R{\mathbb{R}}
\def\W{\mathbb{W}}
\def\Z{\mathbb{Z}}

\def\f{\mib{f}}

\def\v{\mib{v}}
\def\x{\mib{x}}

\def\y{\mib{y}}
\def\vtheta{\mib{\theta}}
\def\0{\mib{0}}

\def\T{{\bf T}}

\def\mC{\mathfrak{C}}
\def\mM{\mathfrak{M}}
\def\mX{\mathfrak{X}}
\def\mY{\mathfrak{Y}}

\def\Det{{\rm Det}}
\def\supp{{\rm supp}\ }

\def\lb{\underline{b}}
\def\ub{\overline{b}}

\theorembodyfont{\itshape}

\newtheorem{thm}{Theorem}[section]
\newtheorem{lem}[thm]{Lemma}

\newtheorem{prop}[thm]{Proposition}

\newtheorem{df}[thm]{Definition}

\newcommand{\mib}[1]{\mbox{\boldmath $#1$}}
\newcommand{\SSC}[1]{\section{#1}\setcounter{equation}{0}}
\newcommand{\qed}{\hbox{\rule[-2pt]{3pt}{6pt}}}

\begin{document}

\begin{Large}
\noindent
{\bf Non-Equilibrium Dynamics of Dyson's Model\\
with an Infinite Number of Particles}
\end{Large}

\vskip 0.5cm

\noindent
{\bf Makoto Katori$^{1}$, Hideki Tanemura$^{2}$} \\

\begin{footnotesize}
\begin{description}
\item{1} \qquad
Department of Physics,
Faculty of Science and Engineering,
Chuo University, \\
Kasuga, Bunkyo-ku, Tokyo 112-8551, Japan.
E-mail: katori@phys.chuo-u.ac.jp 
\item{2} \qquad
Department of Mathematics and Informatics,
Faculty of Science, Chiba University, \\
1-33 Yayoi-cho, Inage-ku, Chiba 263-8522, Japan. 
E-mail: tanemura@math.s.chiba-u.ac.jp
\end{description}
\end{footnotesize}

\begin{center}
(19 June 2009)
\end{center}

\vskip 0.5cm

\noindent{\bf Abstract:}
Dyson's model is a one-dimensional system
of Brownian motions with long-range
repulsive forces acting between any pair
of particles with strength proportional
to the inverse of distances 
with proportionality constant $\beta/2$.
We give sufficient conditions 
for initial configurations so that
Dyson's model with $\beta=2$ and
an infinite number of particles
is well defined in the sense that
any multitime correlation function is
given by a determinant 
with a continuous kernel.
The class of infinite-dimensional configurations
satisfying our conditions is large enough
to study non-equilibrium dynamics.
For example, we obtain the relaxation process 
starting from a configuration, in which
every point of $\Z$ is occupied by
one particle, to the stationary state,
which is the determinantal point process
with the sine kernel.

\pagestyle{plain}

\vskip 0.5cm
\SSC{Introduction}\label{chap: Introduction}

In order to understand the statistics of eigenvalues
of random matrix ensembles as equilibrium distributions
of particle positions in the one-dimensional
Coulomb gas systems with log-potentials,
Dyson introduced stochastic models of particles in $\R$,
which obey the stochastic differential equations (SDEs),
\begin{equation}
dX_j(t)=dB_j(t) 
+ \frac{\beta}{2} \sum_{1 \leq k \leq N, k \not= j}
\frac{dt}{X_j(t)-X_k(t)}, \quad
1 \leq j \leq N, \quad
t \in [0, \infty),
\label{eqn:Dyson}
\end{equation}
where $B_j(t)$'s are independent one-dimensional
standard Brownian motions \cite{Dys62}.
The Gaussian orthogonal ensemble (GOE),
the Gaussian unitary ensemble (GUE),
and the Gaussian symplectic ensemble (GSE)
of random matrices correspond to the SDEs (\ref{eqn:Dyson}) with
$\beta=1, 2$ and 4, respectively \cite{Meh04}.
Spohn \cite{Spo87} has considered an infinite particle system
obtained by taking the $N \to \infty$ limit of 
(\ref{eqn:Dyson}) with $\beta=2$ and called the system
{\it Dyson's model}.
He studied the equilibrium dynamics with respect
to the determinantal (Fermion) point process
$\mu_{\sin}$, in which any spatial correlation
function $\rho_m$ is given by a determinant with
the {\it sine kernel} \cite{Sos00,ST03}
\begin{equation}
K_{\sin}(y-x)= 
\frac{1}{2 \pi} \int_{|k| \leq \pi} dk \,
e^{i k(y-x)}
= \frac{\sin \{ \pi(y-x) \} }{\pi (y-x)},
\quad x, y \in \R,
\label{eqn:sine_kernel0}
\end{equation}
where $i=\sqrt{-1}$.
By the Dirichlet form approach Osada \cite{Osa96}
constructed the infinite particle system
represented by a diffusion process, which has
$\mu_{\sin}$ as a reversible measure.
Recently he proved that this system satisfies
the SDEs (\ref{eqn:Dyson}) with $N= \infty$ \cite{Osa08}.
On the other hand, it was shown by Eynard and Mehta \cite{EM98}
that multitime correlation functions for the process
(\ref{eqn:Dyson}) are generally given by determinants,
if the process starts from $\mu_{N, \sigma^2}^{\rm GUE}$,
the eigenvalue distribution of GUE with variance $\sigma^2$.
Nagao and Forrester \cite{NF98}
evaluated the bulk scaling limit 
$\sigma^2= 2N/\pi^2 \to \infty$ 
and derived the so-called {\it extended sine kernel}
with density 1,
\begin{eqnarray}
{\bf K}_{\sin}(t-s, y-x) &=&
\frac{1}{2 \pi} \int_{|k| \leq \pi} dk \,
e^{k^2(t-s)/2 + i k (y-x)}
- {\bf 1}(s>t) p(s-t, x|y) \nonumber\\
&=& \left\{ \begin{array}{ll} 
\displaystyle{
\int_{0}^{1} du \, e^{\pi^2 u^2 (t-s)/2} 
\cos \{ \pi u (y-x)\} }
& \mbox{if $t>s $} \cr
& \cr
K_{\sin}(y-x)
& \mbox{if $t=s$} \cr
& \cr
\displaystyle{
- \int_{1}^{\infty} du \, 
e^{\pi^2 u^2 (t-s)/2} \cos \{ \pi u (y-x) \} }
& \mbox{if $t<s$},
\end{array} \right.
\label{eqn:sine-kernel}
\end{eqnarray}
$s, t \geq 0, x, y \in \R$,
where ${\bf 1}(\omega)$ is the
indicator function of condition $\omega$,
and $p(t,y|x)$ is the {\it heat kernel}
\begin{equation}
p(t, y|x)=\frac{e^{-(y-x)^2/2t}}{\sqrt{2 \pi t}}
= \frac{1}{2 \pi} \int_{\R} dk \,
e^{-k^2 t/2 + i k (y-x)}, \quad
t > 0.
\label{eqn:p1}
\end{equation}
Since $\lim_{N \to \infty} 
\mu^{\rm GUE}_{N, 2N/\pi^2}=\mu_{\sin}$,
the process, whose multitime correlation functions
are given by determinants with the extended sine 
kernel (\ref{eqn:sine-kernel}),
is expected to be identified with
the infinite-dimensional equilibrium dynamics of Spohn and Osada.
This equivalence is, however, not yet proved. 

Fritz \cite{Fritz87} established
the theory of non-equilibrium dynamics of
infinite particle systems
with a finite-range smooth potential.
Here we study the non-equilibrium dynamics
of infinite-particle Dyson's model with a long-range log-potential, 
in which the force acting each particle is singular
both for short and long distances 
(see (\ref{eqn:Dyson})).

We denote by $\mM$ the space of nonnegative 
integer-valued Radon measures on $\R$,
which is a Polish space with the vague topology:
we say $\xi_n, n \in \N \equiv \{1,2, \dots\}$ 
converges to $\xi$ vaguely, if 
$\lim_{n \to \infty} \int_{\R} \varphi(x) \xi_n(dx)
=\int_{\R} \varphi(x) \xi(dx)$ 
for any $\varphi \in {\rm C}_0(\R)$,
where ${\rm C}_0(\R)$ is the set of all 
continuous real-valued functions with
compact supports.
Any element $\xi$ of $\mM$ can be represented as
$\xi(\cdot) = \sum_{j\in \Lambda}\delta_{x_j}(\cdot)$
with a sequence of points in $\R$, $\x =(x_j)_{j \in \Lambda}$ 
satisfying $\xi(K)=\sharp\{j \in \Lambda: x_j \in K\} < \infty$ 
for any compact subset $K \subset \R$.
The index set $\Lambda$ is $\N$ or a finite set.
We call an element $\xi$ of $\mM$ an unlabeled configuration,
and a sequence $\x$ a labeled configuration.
For $A \subset \R$, we write 
the restriction of $\xi$ on $A$ as 
$(\xi\cap A) (\cdot)=\sum_{j \in \Lambda : x_j \in  A}
\delta_{x_j}(\cdot)$.

As an $\mM$-valued process $(\P, \Xi(t), t\in [0,\infty))$,
we consider the system such that,
for any integer $M \geq 1$,
$f_m \in {\rm C}_{0}(\R), \theta_m \in \R,
1 \leq m \leq M$,
$0 < t_1 < \cdots < t_M < \infty$,
the expectation of
$\exp \Big\{ \sum_{m=1}^{M} \theta_m 
\int_{\R} f_m(x) \Xi(t_m, dx) \Big\}$
can be expanded with $\chi_{m}(x)=e^{\theta_{m} f_{m}(x)}-1,
1 \leq m \leq M$ as
\begin{eqnarray}
{\cal G}^{\xi}[\chi]
&\equiv& \sum_{N_{1} \geq 0} \cdots
\sum_{N_{M} \geq 0}
\prod_{m=1}^{M}\frac{1}{N_m !}
\int_{\R^{N_{1}}} \prod_{j=1}^{N_1} d x_{j}^{(1)}
 \cdots
\int_{\R^{N_{M}}} 
\prod_{j=1}^{N_{M}} d x_{j}^{(M)} \nonumber\\
&& \qquad \times \prod_{m=1}^{M} \prod_{j=1}^{N_{m}} 
\chi_{m} \Big(x_{j}^{(m)} \Big)
\rho\Big( t_{1}, \x^{(1)}_{N_1}; \dots ; t_{M}, \x^{(M)}_{N_M} \Big),
\nonumber
\end{eqnarray}
where $\x^{(m)}_{N_m}$ denotes
$(x^{(m)}_1, \dots, x^{(m)}_{N_m}), 1 \leq m \leq M$.
Here $\rho$'s are locally integrable functions,
which are symmetric in the sense that
$$
\rho(\dots; t_m, \sigma(\x^{(m)}_{N_m}); \dots)
=\rho(\dots; t_m, \x^{(m)}_{N_m}; \dots)
\quad \mbox{with}
\quad \sigma(\x^{(m)}_{N_m})
\equiv (x^{(m)}_{\sigma(1)}, \dots, x^{(m)}_{\sigma(N_m)})
$$
for any permutation $\sigma \in {\cal S}_{N_m},
1 \leq \forall m \leq M$. 
In such a system $\rho( t_{1}, \x^{(1)}_{N_1};
\dots ; t_{M}, \x^{(M)}_{N_M})$ is called
the $(N_1, \dots, N_{M})$-multitime 
correlation function
and ${\cal G}^{\xi}[\chi]$ the generating function 
of multitime correlation functions.
There are no multiple points with probability one for $t > 0$.
Then we assume that there is a function $\mbK(s,x;t,y)$,
which is continuous with respect to
$(x,y) \in \R^2$ for any fixed $(s,t) \in [0, \infty)^2$,
such that
$$
\rho \Big(t_1,\x^{(1)}_{N_1}; \dots;t_M,\x^{(M)}_{N_M} \Big) 
=\det_{
\substack{1 \leq j \leq N_{m}, 1 \leq k \leq N_{n} \\ 1 \leq m, n \leq M}
}
\Bigg[
\mbK(t_m, x_{j}^{(m)}; t_n, x_{k}^{(n)} )
\Bigg]
$$
for any integer $M \geq 1$, 
any sequence $(N_m)_{m=1}^{M}$ of positive integers, and
any time sequence $0 < t_1 < \cdots < t_M < \infty$.
That is, 
the finite dimensional distributions of the process
are determined by the function $\mbK$.
Let $\T=\{t_1, \dots, t_M\}$.
We note that
$\Xi^{\T}=\sum_{t \in \T} \delta_t \otimes \Xi(t)$
is a determinantal (Fermion)
point process on $\T \times \R$
with an operator ${\cal K}$ given by
${\cal K}f(s,x)=\sum_{t \in \T} \int_{\R} dy \,
\mbK(s,x;t, y) f(t,y)$
for $f(t, \cdot) \in {\rm C}_0(\R), t \in \T$.
The process is then said to be 
{\it determinantal with the correlation kernel} $\mbK$.
When ${\cal K}$ is symmetric, Soshnikov \cite{Sos00}
and Shirai and Takahashi \cite{ST03} gave sufficient
conditions for $\mbK$ to be a correlation kernel
of a determinantal point process.
Though such conditions are not known for asymmetric cases,
a variety of processes, which are determinantal
with asymmetric correlation kernels,
have been studied.
As mentioned above
the process $\Xi(t)=\sum_{j=1}^{N} \delta_{X_j(t)}$ with the SDEs
(\ref{eqn:Dyson}) with $\beta=2$ 
starting from its equilibrium measure
$\mu_{N, \sigma^2}^{\rm GUE}$ is an example \cite{EM98}.
The infinite particle system of Nagao and Forrester \cite{NF98}
is also determinantal with 
the extended sine kernel, which is asymmetric
as shown by (\ref{eqn:sine-kernel}).
(For other examples, see, for instance,
\cite{TW04,KT07}.)

In the present paper we first show that,
for any fixed configuration $\xi^N \in \mM$ with
$\xi^N(\R)=N$, 
Dyson's model starting from $\xi^{N}$ is determinantal 
and its correlation kernel
$\mbK^{\xi^{N}}$ is given by using the 
{\it multiple Hermite polynomials} \cite{I05,BK05,IS05} 
(Proposition \ref{Proposition:Finite}).

For $\xi \in \mM$, when 
$\mbK^{\xi \cap [-L, L]}$
converges to a continuous function as $L \to \infty$,
the limit is written as $\mbK^{\xi}$. If 
$\P_{\xi \cap [-L, L]}$ converges
to a probability measure $\P_{\xi}$
on $\mM^{[0, \infty)}$,
which is determinantal with the correlation kernel
$\mbK^{\xi}$, weakly in the sense of finite dimensional
distributions as $L \to \infty$
in the vague topology, we say that the process 
$(\P_\xi, \Xi(t), t\in [0,\infty))$
is {\it well defined with the correlation kernel}
$\mbK^{\xi}$.
(The regularity of the sample paths 
of $\Xi(t)$ will be discussed elsewhere \cite{KT09}.)
In the case $\xi(\R)=\infty$, the process
$(\P_{\xi}, \Xi(t), t\in [0,\infty))$
is Dyson's model with an infinite number of particles.

For $\xi \in \mM$ with $\xi(\{x\}) \leq 1, 
\forall x \in \R$,
we give sufficient conditions so that
the process $(\P_\xi, \Xi(t), t\in [0,\infty))$
is well defined, in which the correlation kernel
is generally expressed using a double
integral with the heat kernels
of an {\it entire function} represented by 
an infinite product (Theorem \ref{Theorem:Infinite1}).
The configuration in which every point of $\Z$
is occupied by one particle,
$\xi^{\Z}(\cdot) \equiv \sum_{\ell \in \Z} \delta_{\ell}(\cdot)$,
satisfies the conditions and we will show that
Dyson's model starting from $\xi^{\Z}$ is determinantal
with the kernel
\begin{eqnarray}
\mbK^{\xi^{\Z}}(s, x; t, y)
&=& {\bf K}_{\sin}(t-s, y-x) \nonumber\\
\label{eqn:Jacobi_kernel}
&& + \frac{1}{2 \pi} \int_{|k| \leq \pi} dk \,
e^{k^2(t-s)/2+i k (y-x)}
\Big\{ \vartheta_3(x- i k s, 2 \pi i s) -1 \Big\}
\\
&=& {\bf K}_{\sin}(t-s, y-x) \nonumber\\
&& + \sum_{\ell \in \Z \setminus \{0\}}
e^{2 \pi i x \ell - 2 \pi^2 s \ell^2}
\int_{0}^{1} du \,
e^{\pi^2 u^2 (t-s)/2}
\cos \Big[ \pi u 
\{ (y-x) - 2 \pi i s \ell \} \Big],
\nonumber
\end{eqnarray}
$s, t \geq 0, x, y \in \R$,
where $\vartheta_3$ is a version of the Jacobi theta function
defined by
\begin{equation}
\vartheta_3(v, \tau) 
= \sum_{\ell \in \Z} e^{2 \pi i v \ell+\pi i \tau \ell^2},
\quad \Im \tau > 0.
\label{eqn:theta3}
\end{equation}
The lattice structure
$\mbK^{\xi^{\Z}}(s, x+n; t, y+n)=\mbK^{\xi^{\Z}}(s,x; t, y),
\forall n \in \Z, s, t \geq 0$ is clear
in (\ref{eqn:Jacobi_kernel}) by the periodicity
of $\vartheta_3$,
$\vartheta_3(v+n, \tau)=\vartheta_3(v,\tau), 
\forall n \in \Z$.
We can prove
\begin{equation}
\lim_{u \to \infty}
\mbK^{\xi^{\Z}}(u+s,x; u+t, y)={\bf K}_{\sin}(t-s, y-x),
\label{eqn:conv}
\end{equation}
which implies that
$\mu_{\sin}$ is an attractor of Dyson's model
and $\xi^{\Z}$ is in its basin.

We are interested in the continuity of the process
with respect to initial configuration.
For Dyson's model with finite particles,
the weak convergence of the processes
$(\P_{\xi^N_n}, \Xi(t), t \in [0, \infty)) \to
(\P_{\xi^N}, \Xi(t), t \in [0, \infty))$ 
as $n \to \infty$ is guaranteed by the 
vague convergence of the initial configurations
$\xi^N_n \to \xi^N$ as $n \to \infty$,
where $\xi^N_n, \xi^N \in \mM$ with
$\xi^N_n(\R)=\xi^N(\R)=N < \infty, n \in \N$.
Based on this continuity, Dyson's model can be
defined for any initial configurations with 
finite particles, which can have multiple points
(see Proposition \ref{Proposition:Finite}).
On the other hand, we have found that,
if $\xi(\R)=\infty$, the weak convergence of processes
in the sense of finite dimensional distributions 
cannot be concluded from the convergence of
initial configurations in the vague topology.
In the present paper we consider a stronger
topology for infinite-particle configurations
(Definition \ref{def:df2}).
We introduce the spaces
$\mY_m^{\kappa}, \kappa \in (1/2,1), m \in \N$
of initial configurations such that 
the convergence of processes is guaranteed by that
of the initial configurations in this new
topology (Theorem \ref{Theorem:Infinite2}).

Note that the union of the spaces
$\mY=\bigcup_{\kappa \in (1/2,1)}
\bigcup_{m \in \N} \mY_m^{\kappa}$
is large enough to carry the Poisson point processes,
Gibbs states with regular conditions, $\mu_{\sin}$,
as well as infinite-particle configurations
with multiple points.
In particular, using the fact $\mu_{\sin}(\mY)=1$ 
and the continuity with respect to the initial configurations, 
we can prove that the process
$({\bf P}_{\sin}, \Xi(t), t \in [0, \infty))$
of Nagao and Forrester, which is determinantal
with the extended sine kernel (\ref{eqn:sine-kernel}),
is Markovian \cite{KT09}.

The paper is organized as follows.
In Section 2 preliminaries and main results
are given.
In Section 3 the definitions of some special functions
used in the present paper are given and their basic
properties are summarized.
Section 4 is devoted to proofs of results.

\SSC{Preliminaries and Main Results}

For $\xi(\cdot)=\sum_{j\in \Lambda}
\delta_{x_j}(\cdot) \in\mM$,
we introduce the following operations;
\begin{description}

\item[(shift)] for $u \in \R$, 
$\tau_u \xi(\cdot) =\displaystyle{\sum_{j \in \Lambda}} 
\delta_{x_j+u}(\cdot)$,

\item[(dilatation)] for $c>0$,
$c \circ \xi(\cdot)=\displaystyle{\sum_{j \in \Lambda} 
\delta_{c x_j}(\cdot)}$,

\item[(square)]
$\displaystyle{
\xi^{\langle 2 \rangle}(\cdot)
=\sum_{j \in \Lambda} \delta_{x_j^2} (\cdot)}$.
\end{description}
We use the convention such that
$$
\prod_{x\in\xi}f(x) =\exp
\left\{\int_\R \xi(dx) \log f(x) \right\}
=\prod_{x \in \supp \xi}f(x)^{\xi(\{x\})}
$$
for $\xi\in \mM$ and a function $f$ on $\R$,
where $\supp \xi = \{x \in \R : \xi(\{x\}) > 0\}$.
For a multivariate symmetric function $g$ we write 
$g((x)_{x \in \xi})$ for $g((x_j)_{j \in \Lambda})$.

For $s,t \in [0,\infty)$, $x,y\in\R$ and 
$\xi^N \in \mM$ with $\xi^N(\R)=N\in\N$,
we set
\begin{eqnarray}
 \mbK^{\xi^N}(s, x; t, y)
&=& \frac{1}{2 \pi i} 
\oint_{\Gamma(\xi^N)} dz \, p(s, x|z)
\int_{\R} dy' \, p(t, -iy|y') \nonumber\\
&& \hskip 2cm \times
\frac{1}{iy'-z}
\prod_{x' \in \xi^{N}}
\left( 1- \frac{iy'-z}{x'-z} \right)
\nonumber\\
&& - {\bf 1}(s > t)p(s-t, x|y),
\label{eqn:KN1a}
\end{eqnarray}
where $\Gamma(\xi^N)$ is a closed contour on the
complex plane $\C$ encircling the points in 
$\supp \xi^N$ on the real line $\R$
once in the positive direction.
\begin{prop}
\label{Proposition:Finite}
Dyson's model $(\P_{\xi^N}, \Xi(t), t\in [0,\infty))$,
starting from any fixed configuration $\xi^N\in \mM$ 
with $\xi^N(\R) = N < \infty$,
is determinantal with the correlation kernel $\mbK^{\xi^N}$
given by (\ref{eqn:KN1a}). 
\label{eqn:THfinite}
\end{prop}
\vskip 0.3cm
We put 
$$
\mM_0= \Big\{ \xi\in\mM : \xi(\{x\})\le 1 \mbox { for any }  x\in\R
\Big\}.
$$ 
Since any element $\xi$ of $\mM_0$ is determined uniquely 
by its support, 
it is identified with a countable subset $\{x_j\}_{j\in\Lambda}$ of $\R$.
For $\xi^{N} \in \mM_0, a \in \C$,
we introduce an entire function of $z \in \C$
$$
\Phi(\xi^N, a, z)= 
\prod_{x \in \xi^{N} \cap \{a\}^{\rm c}}
\left( 1 - \frac{z-a}{x-a} \right),
$$
whose zero set is 
$\supp (\xi^{N} \cap \{a\}^{\rm c})$ 
(see, for instance, \cite{L96}).
Then, if $\xi^N\in \mM_0$,
(\ref{eqn:KN1a}) is written as
\begin{eqnarray}
\mbK^{\xi^N}(s, x; t, y)
&=& \int_{\R} \xi^N(d x') \, p(s, x|x')
\int_{\R} dy' \, p(t, -iy|y')
\Phi(\xi^{N}, x', i y')
\nonumber\\
&&- {\bf 1}(s > t)p(s-t, x|y).
\label{eqn:KN1}
\end{eqnarray}

For $L>0, \alpha>0$ and $\xi\in\mM$ we put
$$
M(\xi, L)=\int_{[-L,L]\setminus\{0\}} \frac{\xi(dx)}{x},
\qquad
M_\alpha(\xi, L)
=\left( \int_{[-L,L]\setminus\{0\}} 
\frac{\xi(dx)}{|x|^\alpha}\right)^{1/\alpha},
$$
and
$$
M(\xi) = \lim_{L\to\infty}M(\xi, L),
\qquad
M_\alpha(\xi)= \lim_{L\to\infty}M_\alpha(\xi, L),
$$
if the limits finitely exist.
We introduce the following conditions:

\vskip 3mm

\noindent ({\bf C.1})
there exists $C_0 > 0$ such that
$|M(\xi)|  < C_0$,

\vskip 3mm

\noindent ({\bf C.2}) (i) 
there exist $\alpha\in (1,2)$ and $C_1>0$ such that
$
M_\alpha(\xi) \le C_1,
$ \\
\noindent (ii) 
there exist $\beta >0$ and $C_2 >0$ such that
$$
M_1(\tau_{-a^2} \xi^{\langle 2 \rangle}) \leq C_2
(|a| \vee 1)^{-\beta}
\quad \forall a \in \supp \xi.
$$

\vskip 3mm

\noindent We denote by $\mX$ the set of configurations $\xi$
satisfying the conditions ({\bf C.1}) and ({\bf C.2}),
and put $\mX_0 = \mX \cap \mM_0$.
For $\xi \in \mX_0$, $a\in\R$ and $z\in\C$ we define
$$
\Phi (\xi,a,z)=\lim_{L\to\infty}
\Phi(\xi \cap [a-L, a+L], a, z).
$$
We note that $|\Phi (\xi,a,z)| < \infty$ 
and $\Phi(\xi, a, \cdot) \not\equiv 0$, 
if $|M(\tau_{-a}\xi)|<\infty$
and $M_2(\tau_{-a} \xi) < \infty$.

\begin{thm}
\label{Theorem:Infinite1}
If $\xi\in \mX_0$, the process $(\P_\xi, \Xi(t), t\in [0,\infty))$
is well defined with the correlation kernel given by
\begin{eqnarray}
\mbK^{\xi}(s, x; t, y)
&=& \int_{\R} \xi(dx') \, p(s, x|x')
\int_{\R} d y' \, p(t, -iy|y')
\Phi (\xi, x',i y')
\nonumber\\
&& - {\bf 1}(s > t) p(s-t, x|y).
\label{eqn:K1}
\end{eqnarray}
\end{thm}
In case $\xi(\R)=\infty$, Theorem \ref{Theorem:Infinite1}
gives Dyson's model with an infinite number of particles
starting from the configuration $\xi \in \mX_0$.
From (\ref{eqn:K1}) it is easy to check that
$$
\mbK^{\xi}(t, x; t, y) \mbK^{\xi}(t, y; t, x) dxdy 
\to \xi(dx){\bf 1}(x=y), 
\quad t \to 0 \quad
\mbox{in the vague topology}.
$$

An interesting and important example is obtained for 
the initial configuration, 
in which every point in $\Z$ is occupied by one particle,
$\xi^{\Z}(\cdot) \equiv \sum_{\ell \in \Z}
\delta_{\ell}(\cdot)$.
In this case 
$\xi^{\Z}(\cdot) \in \mX_0$
and 
we can show that 
the correlation kernel $\mbK^{\xi^{\Z}}$ is given by
(\ref{eqn:Jacobi_kernel}).
The process $({\bf P}_{\sin}, \Xi(t), t\in [0,\infty))$
is reversible with respect to $\mu_{\sin}$.
The result (\ref{eqn:conv}) implies that
the process $(\P_{\xi^{\Z}}, \Xi(u+t), t\in [0,\infty))$
converges to $({\bf P}_{\sin}, \Xi(t), t\in [0,\infty))$,
as $u \to\infty$, weakly in the sense of 
finite dimensional distributions.
In other words, $(\P_{\xi^{\Z}}, \Xi(t), t \in [0, \infty))$
is the {\it relaxation process} from an initial
configuration $\xi^{\Z}$
to the invariant measure $\mu_{\sin}$,
which is determinantal, and this non-equilibrium dynamics
is completely determined via the temporally
inhomogeneous correlation kernel (\ref{eqn:Jacobi_kernel}).
(See Remark in Section 4.3.)

For $\kappa>0$, we put
$$
g^\kappa(x) ={\rm sgn}(x) |x|^{\kappa}, \ x\in\R,
\mbox{ and } 
\eta^{\kappa}(\cdot)=\sum_{\ell\in\Z} 
\delta_{g^\kappa(\ell)}(\cdot).
$$
Since $g^\kappa$ is an odd function, 
$\eta^{\kappa}$ satisfies ({\bf C.1}) 
for any $\kappa>0$.
For any $\kappa > 1/2$ we can show by simple calculation
that $\eta^{\kappa}$ satisfies ({\bf C.2})(i) 
with any $\alpha \in (1/\kappa, 2)$ 
and some $C_1=C_1(\alpha)>0$ depending on $\alpha$,
and does ({\bf C.2})(ii) 
with any $\beta \in (0, 2 \kappa-1)$ 
and some $C_2=C_2(\beta) >0$ depending on $\beta$.
This implies that $\eta^{\kappa}$ is an element of $\mX_0$
in any case $\kappa > 1/2$.
Note that $\eta^{1}=\xi^{\Z}$. 

\vskip 3mm

If there exists $\beta' < (\beta-1) \wedge (\beta/2)$ 
for $\xi \in \mM_0$ such that
$\sharp \{ x \in \xi: 
\xi([x-|x|^{\beta'}, x+|x|^{\beta'}]) \geq 2\} = \infty$, 
then $\xi$ does not satisfy the condition
({\bf C.2}) (ii).
In order to include such initial configurations as well as
those with multiple points
in our study of Dyson's model 
with an infinite number of particles,
we introduce another condition for configurations:

\vskip 3mm

\noindent ({\bf C.3})
there exists $\kappa \in (1/2,1)$ and $m\in\N$ such that
$$
m(\xi,\kappa)\equiv 
\max_{k\in\Z} \xi\bigg( [ g^\kappa(k), g^\kappa(k+1)] \bigg) \le m.
$$

\vskip 3mm

\noindent We denote by $\mY^\kappa_m$ the set of configurations $\xi$
satisfying ({\bf C.1}) and ({\bf C.3}) 
with $\kappa\in (1/2,1)$ and $m\in\N$, and put
$$
\mY = \bigcup_{\kappa\in (1/2,1)}\bigcup_{m\in\N}\mY^\kappa_m.
$$
Noting that the set $\{\xi \in \mM: m(\xi,\kappa)\le m\}$ is 
relatively compact 
for each $\kappa\in (1/2,1)$ and $m\in \N$,
we see that $\mY$ is locally compact.
We introduce the following topology on $\mY$.

\begin{df}
\label{def:df2}
Suppose that $\xi, \xi_n \in \mY, n \in \N$.
We say that $\xi_n$ converges $\Phi$-moderately to $\xi$, if
\begin{equation}
\lim_{n\to\infty} \Phi(\xi_n, i, \cdot)
= \Phi(\xi, i, \cdot)
\mbox{ uniformly on any compact set of $\C$.}
\label{eqn:THC2}
\end{equation}
\end{df}
It is easy to see that (\ref{eqn:THC2}) is satisfied, if
$\xi_n$ converges to $\xi$ vaguely
and the following two conditions hold:
\begin{eqnarray}
&&\lim_{L\to\infty} \sup_{n>0 } 
\Bigg| \int_{[-L,L]^{\rm c}} \frac{\xi_n(dx)}{x} \Bigg|=0,
\label{eqn:THC3}
\\
&&\lim_{L\to\infty} \sup_{n>0 } 
\int_{[-L,L]^{\rm c}} 
\frac{\xi_n^{\langle 2 \rangle}(dx)}{x} =0.
\label{eqn:THC4}
\end{eqnarray}
Note that for any $a\in\R$ and $z\in \C$ 
\begin{equation} 
\lim_{n\to\infty}\Phi(\xi_n,a,z)= \Phi(\xi,a,z), 
\label{con:moderate} 
\end{equation} 
if $\xi_n$ converges $\Phi$-moderately to $\xi$
and $a\notin \supp \xi$.

Then the second theorem of the present paper is the following.

\begin{thm}
\label{Theorem:Infinite2}
{\rm (i)} If $\xi\in\mY$,
$(\P_\xi, \Xi(t), t\in [0,\infty))$
is well defined with a correlation kernel $\mbK^{\xi}$.
In particular, when $\xi \in \mY_{0} \equiv \mY \cap \mM_0$,
$\mbK^{\xi}$ is given by (\ref{eqn:K1}).

\noindent {\rm (ii)} 
Suppose that $\xi, \xi_n \in \mY_m^\kappa, n\in\N$
for some $\kappa\in (1/2,1)$ and $m\in\N$.
If $\xi_n$ converges $\Phi$-moderately to $\xi$, then
the process $(\P_{\xi_n}, \Xi(t), t\in [0,\infty))$
converges to 
the process $(\P_{\xi}, \Xi(t), t\in [0,\infty))$
weakly in the sense of finite dimensional distributions
as $n \to \infty$ in the vague topology.
\end{thm}
In the proof of this theorem given in Section 4.4,
we will give an expression (\ref{eqn:kxi}) to $\mbK^{\xi}$,
which is valid for any $\xi \in \mY$.
There we will use special functions such as
the Hermite polynomials, $H_k, k \in \N_0 \equiv \N \cup \{0\}$,
the complete symmetric functions $h_k, k \in \N_0$,
and the Schur functions $s_{(k|\ell)}, k, \ell \in \N_0$.

\SSC{Special Functions}
\subsection{Multivariate symmetric functions}

For $n\in\N$, let $\lambda=(\lambda_1,\lambda_2,\dots,\lambda_n)$ 
be a partition of length less than or equal to $n$,
and $\delta=(n-1,n-2,\dots,1,0)$.
For $\x=(x_1,x_2,\dots,x_n)$ consider the skew-symmetric polynomial
$$
a_{\lambda+\delta}(\x)
= \det_{1\le j,k \le n}\bigg[x_j^{\lambda_k+n-k}\bigg].
$$
If $\lambda=\emptyset$, it is the Vandermonde determinant,
which is given by the product of difference of variables:
$$
a_{\delta}(\x)
= \det_{1\le j,k \le n}\bigg[x_j^{n-k}\bigg]
=\prod_{1\le j<k \le n}(x_j-x_k).
$$
The {\it Schur function} of the variables $\x=(x_1,x_2,\dots,x_n)$
corresponding to the partition of length $\le n$ is then defined by
$$
s_\lambda(\x)= \frac{a_{\lambda+\delta}(\x)}{a_{\delta}(\x)},
$$
which is a symmetric polynomial of $\x$ \cite{Mac95}.

In the present paper, the following two special cases are considered:

\noindent{\rm (i)} 
When $\lambda=(r)$, $s_{\lambda}(\x)$ is denoted by $h_r(\x)$ and called
the $r$-th {\it complete symmetric function}, which is the sum of all
monomials of total degree $r$ in the variables $\x=(x_1,x_2,\dots,x_n)$.
The generating function for $h_r$ is
$$
H(\x,z)=\sum_{r \in \N_0} h_r(\x)z^r
=\prod_{j=1}^{n} \frac{1}{1-x_jz}
\quad \mbox{ for } \max_{1\le j\le n}|x_jz|<1.
$$

\noindent{\rm (ii)} 
When $\lambda=(k+1,1^\ell)$, $k+\ell+1 \le n$,
we use Frobenius' notation $(k|\ell)$
for the partition,
and consider the Schur function $s_{(k|\ell)}$.
Note that the sum of coefficients 
of the polynomial $s_{(k|\ell)}(\x)$ equals 
\begin{equation}
\label{skl1}
s_{(k|\ell)}(1,\dots,1)= {k+\ell \choose \ell} {n \choose k+\ell+1}.
\end{equation}

Next we consider an infinite sequence of variables:
$\x=(x_j)_{j \in\N}$. 
If $\sum_{j \in \N}x_j <\infty$, and $z$ is a variables such that
$\sup_{j\in\N}|x_j z|<1$, then $|H(\x,z)| <\infty$.
Moreover, if 
$
\sum_{j \in \N}x_j^2 <\infty
$
in addition to the above conditions, we can show
$$
\frac{d^k}{dz^k}\prod_{j \in \N}\frac{1}{1-x_j z}\Bigg|_{z=0}
\le \left.
\frac{d^k}{d z^k}
\exp \left\{
\left|\sum_{j \in\N}x_j \right|z 
+ \sum_{j \in\N}\frac{x_j^2}{1-|x_j z|}z^2
\right\}
\right|_{z=0},
\quad k\in\N,
$$ 
by simple calculation.
It implies
\begin{equation}
\label{complete2}
\sum_{r \in \N_0} |h_r(\x)|z^r 
\le \exp \left\{
\left|\sum_{j \in\N}x_j \right|z 
+ \sum_{j \in\N}\frac{x_j^2}{1-|x_j z|}z^2
\right\},
\end{equation}
and thus the formula
\begin{equation}
\label{complete}
\sum_{r \in \N_0} h_r(\x)z^r = \prod_{j \in\N}\frac{1}{1-x_j z}
\end{equation}
is valid for the infinite sequence of variables $\x=(x_j)_{j\in\N}$. 
Assume that there exist $x_0\in\R$ and $\varepsilon>0$ such that
$\xi\big([x_0-\varepsilon, x_0+\varepsilon]\big)=0$.
We see that for fixed $z\in\C$
\begin{eqnarray}
\Phi(\xi, x,z) &=& \prod_{u \in \xi -\delta_x}
\left( 1- \frac{z-x}{u-x} \right)
\nonumber\\
&=& \prod_{u \in \xi}
\left( 1- \frac{z-x_0}{u-x_0} \right)
\prod_{u \in \xi-\delta_x}\frac{1}{1-(x-x_0)/(u-x_0)}
\nonumber\\
&=& \Phi(\xi, x_0,z)\sum_{r \in \N_0}
h_r\bigg(\Big(\frac{1}{u-x_0}\Big)_{u\in\xi} \bigg)(x-x_0)^r,
\nonumber
\end{eqnarray}
where (\ref{complete}) has been used.
Then $\Phi(\xi, x,z)$ is a smooth function of $x$ on 
$[x_0-\varepsilon, x_0+\varepsilon]$.

\subsection{Multiple Hermite polynomials}

For any $\xi \in \mM$ with $\xi(\R) < \infty$,
the {\it multiple Hermite polynomial of type II},
$P_{\xi}$ is defined as the monic polynomial
of degree $\xi(\R)$ that satisfies for any 
$x \in \supp \xi$
\begin{equation}
\int_{\R} dy \, 
P_{\xi}(y) y^{j} e^{-(y-x)^2/2} = 0,
\quad j=0, \dots, \xi(\{x\})-1.
\label{eqn:Pxi1}
\end{equation}
The {\it multiple Hermite polynomials of type I} 
consist of a set of polynomials
\begin{equation}
\Big\{ A_{\xi}(\, \cdot \, ,x) :
x \in \supp \xi , \quad
{\rm deg} A_{\xi}(\cdot, x)
= \xi(\{x\})-1 \Big\}
\label{eqn:Axi1}
\end{equation}
such that the function
\begin{equation}
Q_{\xi}(y)=\sum_{x \in \supp \xi}
A_{\xi}(y,x) e^{-(y-x)^2/2}
\label{eqn:Qxi1}
\end{equation}
satisfies
\begin{equation}
\int_{\R} dy \,
Q_{\xi}(y) y^{j} = \left\{ \begin{array}{ll}
0, & j=0, \dots, \xi(\R)-2 \cr
1, & j=\xi(\R)-1.
\end{array} \right.
\label{eqn:typeIorth}
\end{equation}
The polynomials $\{A_{\xi}(\cdot , x)\}$
are uniquely determined by 
the degree requirements (\ref{eqn:Axi1})
and the orthogonality relations (\ref{eqn:typeIorth})
\cite{I05}.
The multiple Hermite polynomial of type II, $P_{\xi}$
and the function $Q_{\xi}$ defined by (\ref{eqn:Qxi1})
have the following integration representations \cite{BK05},
\begin{eqnarray}
\label{eqn:Pint}
P_{\xi}(y) &=& 
\int_{\R} d y' \,
\frac{e^{-(y'+iy)^2/2}}{\sqrt{2 \pi}} 
\prod_{x \in \xi}(iy'-x), \cr
\label{eqn:Qint}
Q_{\xi}(y) &=& \frac{1}{2 \pi i}
\oint_{\Gamma(\xi)} d z \,
\frac{e^{-(z-y)^2/2}}{\sqrt{2 \pi}}
\frac{1}{\prod_{x \in \xi} (z-x)}.
\end{eqnarray}

Now we fix $\xi^N \in \mM$ with $\xi^{N}(\R)=N \in \N$.
We write $\xi^N(\cdot) = \sum_{j=1}^{N} \delta_{x_j}(\cdot)$
with a labeled configuration
$\x=(x_j)_{j=1}^{N}$ such that
$x_1 \leq x_2 \leq \cdots \leq x_N$.
Then we define
$$
\xi^{N}_0(\cdot) \equiv 0 \quad
\mbox{and} \quad
\xi^{N}_j(\cdot)= \sum_{k=1}^{j} \delta_{x_k}(\cdot),
\quad 1 \leq j \leq  N.
$$
By definition
$\xi^{N}_{j}(\R)=j, 0 \leq j \leq N$
and $\xi^{N}_j(\{x\}) \leq \xi^{N}_{j+1}(\{x\}),
\forall x \in \R, 0 \leq j \leq N-1$. 
We define
\begin{equation}
H^{(-)}_j(y; \xi^N)=P_{\xi^{N}_j}(y),
\quad
H^{(+)}_j(y; \xi^N)=Q_{\xi^{N}_{j+1}}(y), 
\quad 0 \leq j \leq N-1.
\label{eqn:H+-def}
\end{equation}
By the orthogonality relations (\ref{eqn:Pxi1}),
(\ref{eqn:typeIorth}) and the above definitions, we can prove
the {\it biorthonormality} \cite{BK05}
\begin{equation}
\int_{\R} dy \,
H^{(-)}_j(y; \xi^N) H^{(+)}_k(y; \xi^N)
=\delta_{jk}, \quad
0 \leq j, k \leq N-1.
\label{eqn:orth1}
\end{equation}
For $N \in \N$, let 
$\W_N=\{\x \in \R^N: x_1 < x_2 < \cdots < x_N\}$,
the Weyl chamber of type A$_{N-1}$.

\begin{lem}\label{lemma:detH+}
Let $\y=(y_j)_{j=1}^{N} \in \W_N$.
For any 
$\xi^N (\cdot) = \sum_{j=1}^{N} \delta_{x_j}(\cdot)
\in \mM$ with a labeled configuration
$\x=(x_j)_{j=1}^{N}$ such that 
$x_1 \leq x_2 \leq \cdots \leq x_N$, 
\begin{equation}
\frac{1}{a_{\delta}(\x)}
\det_{1 \leq j, k \leq N} \Big[
e^{-(y_k-x_j)^2/2} \Big]
= (-1)^{N(N-1)/2} (2 \pi)^{N/2}
\det_{1 \leq j, k \leq N}
\Big[ H^{(+)}_{j-1}(y_k; \xi^N) \Big].
\label{eqn:detH+}
\end{equation}
Here when some of the $x_j$'s coincide, 
we interpret the LHS using 
l'H\^opital's rule.
\end{lem}
\noindent{\it Proof.} 
First we assume $\xi^N \in \mM_0$.
Since
$a_{\delta}(\x)=(-1)^{N(N-1)/2}
\prod_{j=2}^{N} \prod_{m=1}^{j-1} (x_j-x_m)$,
by the multilinearity of determinant
\begin{eqnarray}
&& \frac{1}{a_{\delta}(\x)}
\det_{1 \leq j, k \leq N} \Big[
e^{-(y_k-x_j)^2/2} \Big]
\nonumber\\
&=& (-1)^{N(N-1)/2} (2 \pi)^{N/2}
\det_{1 \leq j, k \leq N}
\left[ \frac{e^{-(y_k-x_j)^2/2}}{\sqrt{2 \pi}}
\frac{1}{\prod_{m=1}^{j-1} (x_j-x_m)} \right]
\nonumber\\
&=& (-1)^{N(N-1)/2} (2 \pi)^{N/2}
\det_{1 \leq j, k \leq N}
\left[ \sum_{\ell=1}^{j} 
\frac{e^{-(y_k-x_{\ell})^2/2}}{\sqrt{2 \pi}}
\frac{1}{\prod_{1 \leq m \leq j, m \not=\ell}
(x_{\ell}-x_m)} \right].
\nonumber
\end{eqnarray}
By definition (\ref{eqn:H+-def}) with (\ref{eqn:Qint}),
if $\xi^N \in \mM_0, \xi^N(\R)=N$,
\begin{eqnarray}
\label{eqn:H+j-1A}
H^{(+)}_{j-1}(y_k; \xi^N) &=& 
\frac{1}{2 \pi i}
\oint_{\Gamma(\xi^{N}_{j})} d z \,
\frac{e^{-(y_k-z)^2/2}}{\sqrt{2 \pi}}
\frac{1}{\prod_{x \in \xi^{N}_{j}}(z-x)} \\
&=& 
\frac{1}{2 \pi i}
\oint_{\Gamma(\xi^{N}_{j})} d z \,
\frac{e^{-(y_k-z)^2/2}}{\sqrt{2 \pi}}
\frac{1}{\prod_{\ell=1}^{j} (z-x_{\ell})} 
\nonumber\\
&=& \sum_{\ell=1}^{j} 
\frac{e^{-(y_k-x_{\ell})^2/2}}{\sqrt{2\pi}}
\frac{1}{\prod_{1 \leq m \leq j, m \not= \ell}
(x_{\ell}-x_m)}, \quad 1 \leq j \leq N.
\nonumber
\end{eqnarray}
Then (\ref{eqn:detH+}) is proved for $\xi^N \in \mM_0$.
When some of the $x_j$'s coincide,
the LHS of (\ref{eqn:detH+}) is interpreted
using l'H\^opital's rule and
in the RHS of (\ref{eqn:detH+})
$H^{(+)}_{j-1}(y_k; \xi^N)$ should be given by
(\ref{eqn:H+j-1A}).
Then (\ref{eqn:detH+}) is valid for any
$\xi^N \in \mM, \xi^N(\R)=N$. \qed
\vskip 0.5cm
\begin{lem}\label{lemma:H+-p}
Let $N \in \N, \xi^N \in \mM$ with $\xi^N(\R) =N$.
For $0 \leq s \leq t, x, y \in \R, 0 \leq j \leq N-1$,
\begin{eqnarray}
&& \int_{\R} dy \,
H^{(-)}_j \left( \frac{y}{\sqrt{t}}; 
\frac{1}{\sqrt{t}} \circ \xi^N \right) p(t-s, y|x)
=\left( \frac{s}{t} \right)^{j/2}
H^{(-)}_j \left( \frac{x}{\sqrt{s}};
\frac{1}{\sqrt{s}} \circ \xi^N \right), \nonumber\\
\label{eqn:Hminusp}
\\
&& \int_{\R} dx \,
p(t-s, y|x)
H^{(+)}_j \left( \frac{x}{\sqrt{s}}; 
\frac{1}{\sqrt{s}} \circ \xi^N \right)
= \left( \frac{s}{t} \right)^{(j+1)/2}
H^{(+)}_j \left( \frac{y}{\sqrt{t}};
\frac{1}{\sqrt{t}} \circ \xi^N \right),
\nonumber\\
\label{eqn:Hplusp}
\end{eqnarray}
where $p$ is the heat kernel (\ref{eqn:p1}).
\end{lem}
\noindent{\it Proof.} 
Consider the integral 
\begin{eqnarray}
&& \int_{\R} dy \,
H^{(-)}_j \left( \frac{y}{\sqrt{t}};
\frac{1}{\sqrt{t}} \circ \xi^N \right) p(t-s, y|x)
\nonumber\\
&=& \frac{1}{\sqrt{2 \pi (t-s)}}
\frac{1}{\sqrt{2 \pi}}
\int_{\R} d y' \prod_{x \in \xi^{N}_j}
\left( i y' - \frac{x}{\sqrt{t}} \right)
\int_{\R} dy \,
e^{-(y-x)^2/\{2(t-s)\}
-(y'+iy/\sqrt{t})^2/2}
\nonumber\\
&=& \sqrt{\frac{t}{s}} \frac{1}{\sqrt{2\pi}}
\int_{\R} d y' \,
\prod_{x \in \xi^{N}_j} \left( iy' - \frac{x}{\sqrt{t}} \right)
e^{-t(y'+ix/\sqrt{t})^2/(2s)}.
\nonumber
\end{eqnarray}
Change the integral variable $y' \to y' \sqrt{t/s}$
to obtain the equality (\ref{eqn:Hminusp}).
Similar calculation gives (\ref{eqn:Hplusp}). \qed

When $\xi^N(\cdot)=N \delta_0(\cdot)$,
\begin{eqnarray}
H^{(-)}_{j}(y; N \delta_0) &=& 2^{-j/2} H_j(y/\sqrt{2}),
\nonumber\\
H^{(+)}_{j}(y; N \delta_0) &=&
\frac{2^{-j/2}}{j! \sqrt{2 \pi}}
H_j(y/\sqrt{2}) e^{-y^2/2},
\quad 0 \leq j \leq N-1,
\nonumber
\end{eqnarray}
where $H_j(x)$ is the Hermite polynomial of degree $j$,
\begin{eqnarray}
H_j(x) &=& j ! \sum_{k=0}^{[j/2]}
(-1)^k \frac{(2x)^{j-2k}}{k ! (j-2k)!}
\nonumber\\
&=& 2^{j/2}
\int_{\R} dy \,
\frac{e^{-y^2/2}}{\sqrt{2 \pi}} (iy + \sqrt{2} x)^j
\nonumber\\
\label{eqn:Hermite2}
&=& \frac{j!}{2 \pi i} \oint_{\Gamma(\delta_0)}
d z \,
\frac{e^{2 z x-z^2}}{z^{j+1}}.
\end{eqnarray}
The last expression (\ref{eqn:Hermite2}) implies that
the generating function of the Hermite polynomials
is given by
\begin{equation}
e^{2 z x - z^2}
= \sum_{j\in\N_0} \frac{z^{j}}{j !} H_j(x).
\label{eqn:Hermite3}
\end{equation}

\SSC{Proofs of Results}
\subsection{Proof of Proposition \ref{Proposition:Finite}}%

For $\x, \y \in \W_N$ and $t >0$, consider the
{\it Karlin-McGregor determinant} of
the heat kernel (\ref{eqn:p1}) \cite{KM59}
$$
f_N(t, \y|\x)=\det_{1 \leq j, k \leq N}
\Big[ p(t, y_j|x_k) \Big].
$$
If $\xi^N \in \mM_0$ with $\xi^N(\R)=N \in \N$,
$\xi^N$ can be identified with a set $\x \in \W_N$.
For any $M \geq 1$ and any time sequence 
$0 < t_1 < \cdots < t_M < \infty$,
the {\it multitime probability density} of Dyson's model
is given by \cite{Gra99,KT07}
\begin{eqnarray}
&& p^{\xi^N} \Big(t_1, \xi^{(1)}; \dots ;
t_M, \xi^{(M)} \Big) \nonumber\\
&& \qquad 
= a_{\delta}(\x^{(M)}) 
\prod_{m=1}^{M-1} f_N (t_{m+1}-t_m, \x^{(m+1)}|\x^{(m)})
f_N(t_1, \x^{(1)}|\x) \frac{1}{a_{\delta}(\x)},
\nonumber
\end{eqnarray}
where $\xi^{(m)}(\cdot)=\sum_{j=1}^N \delta_{x_j^{(m)}}(\cdot)$,
$1 \leq m \leq M$.

Define 
\begin{eqnarray}
\phi_{j}^{(-)}(t,x; \xi^N)
&\equiv& t^{j/2}
H_{j}^{(-)}\left( \frac{x}{\sqrt{t}}; 
\frac{1}{\sqrt{t}} \circ \xi^N \right),
\nonumber\\
\phi_{j}^{(+)}(t, x; \xi^N)
&\equiv& 
t^{-(j+1)/2}
H_{j}^{(+)} \left( \frac{x}{\sqrt{t}}; 
\frac{1}{\sqrt{t}} \circ \xi^N \right),
\nonumber
\end{eqnarray}
$0 \leq j \leq N-1, t>0, x \in \R$. 
From the biorthonormality (\ref{eqn:orth1})
of the multiple Hermite polynomials
and Lemma \ref{lemma:H+-p}, the following relations
are derived.

\begin{lem}\label{lemma:phi+-p}
For $\xi^N \in \mM$ with $\xi^N(\R)=N \in \N$,
$0 \leq t_1 \leq t_2$,
\begin{eqnarray}
&& \int_{\R} dx_2 \,
\phi_{j}^{(-)}(t_2, x_2; \xi^N) p(t_2-t_1, x_2|x_1)
=\phi_{j}^{(-)}(t_1, x_1; \xi^N), \quad
0 \leq j \leq N-1, \quad \nonumber\\
&& \int_{\R} dx_1 \,
p(t_2-t_1, x_2 |x_1)
\phi_{j}^{(+)}(t_1, x_1; \xi^N)
=\phi_{j}^{(+)}(t_2, x_2; \xi^N), \quad
0 \leq j \leq N-1, \quad \nonumber\\
&& \int_{\R} dx_1 
\int_{\R} dx_2 \,
\phi_{j}^{(-)}(t_2, x_2; \xi^N)
p(t_2-t_1, x_2|x_1)
\phi_{k}^{(+)}(t_1, x_1; \xi^N)
=\delta_{jk}, 
\nonumber\\
&& \hskip 10cm
0 \leq j, k \leq N-1. \nonumber
\end{eqnarray}
\end{lem}
\vskip 0.5cm
Put
$$
\mu^{(\pm)}(t, \x; \xi^N)
=\det_{1 \leq j, k \leq N}
\Big[ \phi^{(\pm)}_{j-1}(t, x_k; \xi^N) \Big].
$$
Since $H^{(-)}_j$ is a monic polynomial of degree $j$,
$\mu^{(-)}(t, \x; \xi^N)
=(-1)^{N(N-1)/2} a_{\delta}(\x)$.
By Lemma \ref{lemma:detH+}, $f_N(t_1, \x^{(1)}|\x)/a_{\delta}(\x)$
will be replaced by 
$(-1)^{N(N-1)/2} \mu^{(+)}(t_1, \x^{(1)}; \xi^N)$
to extend the expression to the case $\xi^N \in \mM$.
Then the multitime probability density of Dyson's model 
is expressed as
\begin{eqnarray}
&& p^{\xi^N} \Big(t_1, \xi^{(1)}; \dots ;
t_M, \xi^{(M)} \Big) \nonumber\\
&& 
=\mu^{(-)}(t_M, \x^{(M)}; \xi^N)
\prod_{m=1}^{M-1} f_N(t_{m+1}-t_m; \x^{(m+1)}|\x^{(m)})
\mu^{(+)}(t_1, \x^{(1)}; \xi^N)
\label{eqn:pxiN2}
\end{eqnarray}
for $\xi^N \in \mM$ with $\xi^N(\R)=N \in \N$.
For $\x=(x_1, \dots, x_N)$ with
$\xi(\cdot)=\sum_{j=1}^{N} \delta_{x_j} (\cdot)$
and $N' \in \{ 1,2, \dots, N\}$, we put
$\x_{N'}=(x_1, \dots, x_{N'})$.
For a sequence $(N_m)_{m=1}^{M}$ of positive integers 
less than or equal to $N$,
we obtain the 
$(N_1, \dots, N_{M})$-{\it multitime correlation function} by
\begin{eqnarray}
&& \rho_N^{\xi^N} \Big(t_{1}, \x^{(1)}_{N_1} ; 
\dots; t_M, \x^{(M)}_{N_M} \Big) 
\nonumber\\
&&=
\int_{\prod_{m=1}^{M} \R^{N-N_{m}}}
p^{\xi^N} \Big(t_1, \xi^{(1)}; 
\dots; t_M, \xi^{(M)} \Big)
\prod_{m=1}^{M}
\frac{1}{(N-N_{m})!}
\prod_{j=N_{m}+1}^{N} dx_{j}^{(m)}. \qquad
\label{def:corr}
\end{eqnarray}

For $\f=(f_{1}, \cdots, f_{M}) \in {\rm C}_{0}(\R)^{M}$,
and $\vtheta=(\theta_{1}, \cdots, \theta_{M})
\in \R^{M}$,
the generating function for multitime correlation functions
is given as
\begin{eqnarray}
{\cal G}^{\xi^N}[\chi]
&=& \E_{\xi^N} \left[
\exp \left\{ \sum_{m=1}^{M}
\theta_{m} \sum_{j_{m}=1}^{N} 
f_{m}(X_{j_{m}}(t_{m})) \right\} \right]
\nonumber\\
&=& \sum_{N_{1}=0}^{N} \cdots
\sum_{N_{M}=0}^{N}
\prod_{m=1}^{M}\frac{1}{N_m !}
\int_{\R^{N_{1}}} \prod_{j=1}^{N_1} d x_{j}^{(1)}
 \cdots
\int_{\R^{N_{M}}} 
\prod_{j=1}^{N_{M}} d x_{j}^{(M)} \nonumber\\
&& \times \prod_{m=1}^{M} \prod_{j=1}^{N_{m}} 
\chi_{m} \Big(x_{j}^{(m)} \Big)
\rho^{\xi^N} \Big( t_{1}, \x^{(1)}_{N_1};
\dots ; t_{M}, \x^{(M)}_{N_M} \Big),
\nonumber
\end{eqnarray}
where
$$
\chi_{m}(x)=e^{\theta_{m} f_{m}(x)}-1, 
\qquad 1 \leq m \leq M.
$$
By the argument given in Section 4.2 in \cite{KT07},
the expression (\ref{eqn:pxiN2}) with
Lemma \ref{lemma:phi+-p} leads to
the Fredholm determinantal expression
for the generating function,
$$
{\cal G}^{\xi^N}[\chi]=
\Det \Big[ \delta_{m n} \delta (x-y)
+\widetilde{S}^{m,n}(x,y; \xi^N)
\chi_n (y) \Big],
$$
where
$$
\widetilde{S}^{m,n}(x,y; \xi^N)
=S^{m,n}(x,y; \xi^N)
-{\bf 1}(m > n) p(t_m-t_n, x|y)
$$
with
\begin{eqnarray}
&& S^{m,n}(x,y; \xi^N)
= \sum_{j=0}^{N-1} 
\phi_{j}^{(+)}(t_m, x; \xi^N)
\phi_{j}^{(-)}(t_n, y; \xi^N)
\nonumber\\
&& \qquad =\frac{1}{\sqrt{t_m}}
\sum_{j=0}^{N-1}
\left( \frac{t_n}{t_m} \right)^{j/2}
H_j^{(+)}\left( \frac{x}{\sqrt{t_m}}; 
\frac{1}{\sqrt{t_m}} \circ \xi^N \right)
H_j^{(-)} \left( \frac{y}{\sqrt{t_n}};
\frac{1}{\sqrt{t_n}} \circ \xi^N \right).
\nonumber
\end{eqnarray}
Here the Fredholm determinant is expanded as
\begin{eqnarray}
&& \Det \Big[
\delta_{m n} \delta(x-y)
+\widetilde{S}^{m,n}(x,y; \xi^N) \chi_{n}(y) \Big]
\nonumber\\
&=& \sum_{N_{1}=0}^{N} \cdots
\sum_{N_{M}=0}^{N}
\prod_{m=1}^{M}\frac{1}{N_m !}
\int_{\R^{N_{1}}} \prod_{j=1}^{N_1} d x_{j}^{(1)}
 \cdots
\int_{\R^{N_{M}}} 
\prod_{j=1}^{N_{M}} d x_{j}^{(M)} \nonumber\\
&& \quad \quad
\times \prod_{m=1}^{M} \prod_{j=1}^{N_{m}} 
\chi_{m} \Big(x_{j}^{(m)} \Big) 
\det_{
\substack{1 \leq j \leq N_{m}, 1 \leq k \leq N_{n} \\ 1 \leq m, n \leq M}
}
 \Bigg[
\widetilde{S}^{m,n}(x^{(m)}_j, x^{(n)}_k; \xi^N)
\Bigg].
\nonumber
\end{eqnarray}

\noindent{\it Proof of Proposition \ref{Proposition:Finite}.} 
Inserting the integral formulas
for $H_j^{(\pm)}$, 
the kernel $S^{m,n}$
is written as 
\begin{eqnarray}
S^{m,n}(x,y; \xi^N) 
&=& \frac{1}{\sqrt{t_m}} 
\frac{1}{2\pi i} \oint_{\Gamma(t_m^{-1/2} \circ \xi^N)} dz \,
\frac{e^{-(z-x/\sqrt{t_m})^2/2}}
{\sqrt{2 \pi}}
\int_{\R} dy' \,
\frac{e^{-(y'+iy/\sqrt{t_n})^2/2}}{\sqrt{2 \pi}}
\nonumber\\
&& \qquad \qquad \times
\sum_{k=0}^{N-1} \left( \frac{t_n}{t_m} \right)^{k/2}
\frac{\prod_{\ell=1}^{k}(iy'-x_{\ell}/\sqrt{t_n})}
{\prod_{\ell=1}^{k+1}(z-x_{\ell}/\sqrt{t_m})}
\nonumber\\
&=& \frac{1}{2\pi i} \oint_{\Gamma(t_m^{-1/2} \circ \xi^N)} dz \,
\frac{e^{-(z-x/\sqrt{t_m})^2/2}}
{\sqrt{2 \pi}}
\int_{\R} dy' \,
\frac{e^{-(y'+iy/\sqrt{t_n})^2/2}}{\sqrt{2 \pi}}
\nonumber\\
&& \qquad \qquad \times
\sum_{k=0}^{N-1}
\frac{\prod_{\ell=1}^{k}(i\sqrt{t_n} y'-x_{\ell})}
{\prod_{\ell=1}^{k+1}(\sqrt{t_m} z-x_{\ell})}.
\nonumber
\end{eqnarray}
For $z_1, z_2 \in \C$
with $z_1 \notin \{x_1, \dots, x_N\}$, the following identity
holds,
\begin{eqnarray}
&& \sum_{k=0}^{N-1} 
\frac{\prod_{\ell=1}^{k} (z_2-x_{\ell})}
{\prod_{\ell=1}^{k+1} (z_1-x_{\ell})}
\nonumber\\
&=& \frac{1}{z_1-x_1}
+\frac{z_2-x_1}{(z_1-x_1)(z_1-x_2)}+ \cdots
+ \frac{(z_2-x_1)(z_2-x_2)\cdots (z_2-x_{N-1})}
{(z_1-x_1)(z_1-x_2)\cdots (z_1-x_N)}
\nonumber\\
&=& \left( \prod_{\ell=1}^{N}
\frac{z_2-x_{\ell}}{z_1-x_{\ell}}
-1 \right)
\frac{1}{z_2-z_1}.
\nonumber
\end{eqnarray}
By this identity, we have
\begin{eqnarray}
S^{m,n}(x,y;\xi^N) 
&=& \frac{1}{2\pi i} \oint_{\Gamma(t_m^{-1/2} \circ \xi^N)} dz \,
\frac{e^{-(z-x/\sqrt{t_m})^2/2}}
{\sqrt{2 \pi}}
\int_{\R} dy' \,
\frac{e^{-(y'+iy/\sqrt{t_n})^2/2}}{\sqrt{2 \pi}}
\nonumber\\
&& \qquad \times
\left( \prod_{\ell=1}^{N}
\frac{i \sqrt{t_n} y' - x_{\ell}}{\sqrt{t_m} z - x_{\ell}}
-1 \right)
\frac{1}{i\sqrt{t_n}y' -\sqrt{t_m} z}.
\nonumber
\end{eqnarray}
Note that
\begin{eqnarray}
&& \frac{1}{2\pi i} \oint_{\Gamma(t_m^{-1/2} \circ \xi^N)} dz \,
\frac{e^{-(z-x/\sqrt{t_m})^2/2}}
{\sqrt{2 \pi}}
\int_{\R} dy' \,
\frac{e^{-(y'+iy/\sqrt{t_n})^2/2}}{\sqrt{2 \pi}}
\frac{1}
{i \sqrt{t_n}y' - \sqrt{t_m} z}
\nonumber\\
&=& 
\frac{1}{2\pi i} \oint_{\Gamma(t_m^{-1/2} \circ \xi^N)} dz \,
\frac{e^{-(z-x/\sqrt{t_m})^2/2}}
{\sqrt{2 \pi}}
\int_{\R} dy' \,
\frac{e^{-(y'+iy/\sqrt{t_n})^2/2}}{\sqrt{2 \pi}}
\frac{1}{i \sqrt{t_n} y'}
\sum_{j \in \N_0}
\left( \sqrt{\frac{t_m}{t_n}} \frac{z}{iy'} \right)^j 
\nonumber\\
&=& 0.
\nonumber
\end{eqnarray}
By changing the integral variables appropriately,
we find that $\widetilde{S}^{m,n}(x,y;\xi^N)$ is
equal to (\ref{eqn:KN1a})
with $s=t_m, t=t_n$. This completes the
proof. \qed
\vskip 0.5cm

\subsection{Proof of Theorem \ref{Theorem:Infinite1}}

In this subsection we give a proof of Theorem 
\ref{Theorem:Infinite1}.
First we prove some lemmas.

\begin{lem}\label{lemma:4_2_1}
If $M_\alpha(\xi) <\infty$ for some $\alpha \in (1,2)$, then
$$
\alpha \sum_{L \in \N}
\frac{M_1(\xi,L)^{\alpha/(\alpha-1)}}{L(L+1)^{\alpha}}
\le M_\alpha(\xi)^{\alpha^2/(\alpha-1)}.
$$
\end{lem}

\noindent {\it Proof.}
By H\"older's inequality we have
$$
M_1(\xi, L)= 
\int_{0 < |x| \leq L} 
\frac{\xi(dx)}{|x|}
\le M_\alpha(\xi)\xi\big([-L,L]\setminus\{0\}\big)^{(\alpha-1)/\alpha}.
$$
On the other hand
\begin{eqnarray}
M_\alpha (\xi)^\alpha
&=& \sum_{L \in \N}
\int_{L-1 < |x| \leq L} 
\frac{\xi(dx)}{|x|^\alpha}
\nonumber\\
&\ge& \sum_{L \in \N} L^{-\alpha}
\bigg\{\xi\big([-L,L]\setminus\{0\}\big)
-\xi\big([-L+1,L-1]\setminus\{0\}\big)\bigg\}
\nonumber\\
&=& \sum_{L \in \N} \Big\{L^{-\alpha}-(L+1)^{-\alpha} \Big\}
\xi\big([-L,L]\setminus\{0\}\big)
\nonumber\\
&\ge& \alpha \sum_{L \in \N}
\frac{\xi\big([-L,L]\setminus\{0\}\big)}
{L (L+1)^{\alpha}}.
\nonumber
\end{eqnarray}
From the above inequalities we have
$$
M_\alpha (\xi)^\alpha 
\ge \alpha \sum_{L \in \N}
\frac{1}{L(L+1)^{\alpha}}
\left( \frac{M_1(\xi,L)}{M_\alpha(\xi)} \right)^{\alpha/(\alpha-1)}.
$$
Lemma \ref{lemma:4_2_1} is derived from this inequality, since 
$\alpha+\alpha/(\alpha-1)
=\alpha^2/(\alpha-1)$.
\qed

\begin{lem}\label{lemma:4_2_2}
Let $\alpha\in (1,2)$ and $\delta> \alpha-1$.
Suppose that $M_\alpha(\xi)<\infty$ and put 
$L_0=L_0(\alpha, \delta, \xi)= 
(2M_\alpha(\xi))^{\alpha/(\delta-\alpha+1)}$.
Then
$$
M_1 (\xi,L) \le L^\delta, \quad L\ge L_0.
$$
\end{lem}

\noindent {\it Proof.} 
Suppose that $L_1\in \N$ satisfies $M_1(\xi,L_1)>L_1^\delta$.
Then
\begin{eqnarray}
&&\alpha \sum_{L \in \N}
\frac{M_1(\xi,L)^{\alpha/(\alpha-1)}}{L(L+1)^{\alpha}}
> \alpha \sum_{L=L_1}^\infty
\frac{L_1^{\alpha\delta/(\alpha-1)}}{L(L+1)^{\alpha}}
\nonumber\\
&&\qquad > \alpha L_1^{\alpha\delta/(\alpha-1)}
\int_{L_1+1}^\infty dy \ y^{-(\alpha+1)}
\nonumber\\
&&\qquad = L_1^{\alpha\delta/(\alpha-1)}(L_1+1)^{-\alpha}
= \left(\frac{L_1}{L_1+1}\right)^\alpha
L_1^{\alpha(\delta-\alpha+1)/(\alpha-1)}.
\nonumber
\end{eqnarray}
From Lemma \ref{lemma:4_2_1} we have
$$
\left(\frac{L_1}{L_1+1}\right)^\alpha
L_1^{\alpha(\delta-\alpha+1)/(\alpha-1)}
\le M_\alpha(\xi)^{\alpha^2/(\alpha-1)}.
$$
Hence
$$
L_1 < \left(\frac{L_1+1}{L_1}\right)^{(\alpha-1)/(\delta-\alpha+1)}
M_\alpha(\xi)^{\alpha/(\delta-\alpha+1)}
< (2M_\alpha(\xi))^{\alpha/(\delta-\alpha+1)}.
$$
This completes the proof. \qed

\vskip 3mm

The following lemma will play an important role
in the proof of Theorem \ref{Theorem:Infinite1}.
\begin{lem}
\label{lemma:4_2_5}
For any $\xi\in\mX_0$, there exist
$C_3= C_3(\alpha, \beta, C_0, C_1, C_2) >0$ 
and $\theta \in (\alpha \vee (2-\beta), 2)$
such that
$$
|\Phi (\xi,a, iy)|
\le \exp \Big[ C_3 \Big\{(|y|^{\theta} \vee 1)
+ (|a|^{\theta} \vee 1) \Big\} \Big]
\quad \forall y \in \R, \; \forall a\in \supp \xi.
$$
\end{lem}

\vskip 3mm
\noindent {\it Proof.}
First we estimate the entire function $\Phi(\xi,a,z), z\in \C$,
in the case that $a=0 \in \supp \xi$.
In case $2|z|<|x|$,
by using the expansion
$$
\log \left(1+\frac{z}{x}\right)
= \sum_{k \in \N} \frac{(-1)^{k-1}}{k}
\left(\frac{z}{x}\right)^k,
$$
we have
\begin{eqnarray}
&&\int_{2|z|<|x|}\xi(dx) \log\left(1+\frac{z}{x}\right)
= \int_{2|z|<|x|}\xi(dx) \sum_{k \in \N}
\frac{(-1)^{k-1}}{k}\left(\frac{z}{x}\right)^k
\nonumber\\
&&\qquad= \int_{2|z|<|x|}\xi(dx) \frac{z}{x}
+ \int_{2|z|<|x|}\xi(dx)
\left(\frac{z}{x}\right)^2
\sum_{k=2}^\infty 
\frac{(-1)^{k-1}}{k}\left(\frac{z}{x}\right)^{k-2}.
\nonumber
\end{eqnarray}
Since
$$
\left|\int_{2|z|<|x|}\xi(dx) \frac{z}{x}\right|
\le |M(\xi)||z| + M_1(\xi,2|z|)|z|,
$$
and
\begin{eqnarray}
&&\left|\int_{2|z|<|x|}\xi(dx) \left(\frac{z}{x}\right)^2
\sum_{k=2}^\infty 
\frac{(-1)^{k-1}}{k}\left(\frac{z}{x}\right)^{k-2}\right|
\nonumber\\
&&\qquad 
\le \int_{2|z|<|x|}\xi(dx) \frac{|z|^2}{|x|^2}
\frac{1}{2}\sum_{k=2}^\infty 2^{2-k}
= \int_{2|z|<|x|}\xi(dx)\frac{|z|^2}{|x|^2}
\nonumber\\
&&\qquad
\le M_\alpha(\xi)^\alpha |z|^\alpha,
\nonumber
\end{eqnarray}
we have
\begin{equation}
\prod_{x\in\xi}
\left\{ 1 + {\bf 1}(|x|>2|z|)\frac{z}{x}\right\}
\le \exp \Big\{
|M(\xi)||z|+M_1(\xi,2|z|)|z|
+M_\alpha(\xi)^\alpha|z|^\alpha
\Big\}.
\label{estimate1}
\end{equation}
On the other hand we have
\begin{eqnarray}
&&\prod_{x\in\xi}
\left\{ 1 + {\bf 1}(0<|x|\le2|z|)\frac{z}{x}\right\}
\nonumber\\
&&\qquad \le 
\exp \left\{\int_{[-2|z|,2|z|]\setminus\{0\}} 
\frac{\xi(dx)}{|x|}|z| \right\}
= \exp \Big\{M_1(\xi,2|z|)|z| \Big\}.
\label{estimate2}
\end{eqnarray}
Combining the above two inequalities
(\ref{estimate1}) and (\ref{estimate2}), we obtain 
$$
\prod_{x \in \xi \cap \{0\}^{\rm c}}
\left( 1 + \frac{z}{x}\right)
\le \exp \Big\{
|M(\xi)||z|+2M_1(\xi,2|z|)|z|
+M_\alpha(\xi)^\alpha|z|^\alpha \Big\}.
$$
By the conditions ({\bf C.1}), ({\bf C.2})(i) and
Lemma \ref{lemma:4_2_2}, we have
$$
|M(\xi)||z|+2M_1(\xi,2|z|)|z|
+M_\alpha(\xi)^\alpha|z|^\alpha
\le C_0|z| + 4 |z|^{1+\delta}+C_1 |z|^{\alpha}
$$
for $|z| \le L_0(\alpha, \delta, C_1)$ with 
$\delta > \alpha-1$.
Hence, if we can take $\theta > \alpha$, then
\begin{equation}
|\Phi (\xi,0, z)|
\le \exp \Big[ C (|z|^{\theta} \vee 1) \Big]
\label{eqn:PhiZ0}
\end{equation}
with a positive constant $C$, 
which depends on only $\alpha, \beta, C_0$
and $C_1$.

Next we estimate the case that $a\not= 0$
by using the following equations:
$$
\Phi(\xi, a, z) = \Phi(\xi, 0, z)
\Phi(\xi \cap \{0\}^{\rm c}, a, 0)
\left(\frac{z}{a}\right)^{\xi(\{0 \})}
\frac{a}{a-z}
$$
and
$$
\Phi(\xi \cap \{0\}^{\rm c},a,0)
=\Phi(\xi \cap \{-a\}^{\rm c},0,-a)
\Phi(\xi^{\langle 2 \rangle} \cap \{0\}^{\rm c}, a^2,0)
2^{1-\xi(\{-a\})},
$$
where $a \in \supp \xi$.
From the above estimate (\ref{eqn:PhiZ0}) we have
$$
|\Phi(\xi \cap \{-a\}^{\rm c},0,-a)|
\leq \exp \Big[ C (|a|^{\theta} \vee 1) \Big].
$$
Since $|(iy/a)^{\xi(\{0\})} a/(a-iy)| \leq 1$, 
it is enough to show the estimate
$$
\Phi(\xi^{\langle 2 \rangle} \cap \{0\}^{\rm c},a^2, 0) 
\leq \exp \Big\{ 3 C_2 |a|^{2-\beta} \Big\}
$$
for proving this lemma.
We have
\begin{eqnarray}
&& \left| \int_{2a^2 < |x-a^2|} \xi^{\langle 2 \rangle}(dx)
\log \left( 1+\frac{a^2}{x-a^2} \right) \right|
\nonumber\\
&& \qquad
= \left|\int_{2a^2 <|x-a^2|} \xi^{\langle 2 \rangle}(dx)
\frac{a^2}{x-a^2} 
\sum_{k \in \N} \frac{(-1)^{k}}{k}
\left( \frac{a^2}{x-a^2} \right)^{k-1} \right|
\nonumber\\
&& \qquad \leq 2 \int_{2 a^2 <|x-a^2|} \xi^{\langle 2 \rangle}(dx)
\left| \frac{a^2}{x-a^2} \right|
\leq 2 M_1(\tau_{-a^2} \xi^{\langle 2 \rangle}) a^2.
\nonumber
\end{eqnarray}
On the other hand we see
\begin{eqnarray}
&& \prod_{x \in \xi^{\langle 2 \rangle}} \left\{
1+{\bf 1}(0 < |x-a^2| < 2 a^2) \frac{a^2}{x-a^2} \right\}
\nonumber\\
&& \qquad \leq
\exp \left\{ \int_{[-2a^2, 2 a^2]\setminus\{0\}}
\frac{(\tau_{-a^2} \xi^{\langle 2 \rangle})(dx)}{|x|} a^2 \right\}
=\exp \Big\{ M_1(\tau_{-a^2} \xi^{\langle 2 \rangle}, 2 a^2) a^2 \Big\}.
\nonumber
\end{eqnarray}
Then
\begin{eqnarray}
\Phi(\xi^{\langle 2 \rangle} \cap \{0\}^{\rm c},
a^2, 0) 
&=& 
\prod_{x \in \xi^{\langle 2 \rangle}\cap \{0, a^2\}^{\rm c}} 
\left(1+\frac{a^2}{x-a^2} \right)
\nonumber\\
&\leq& \exp \Big\{ 3M_1(\tau_{-a^2} \xi^{\langle 2 \rangle}) a^2 \Big\}
=\exp \Big\{ 3 C_2 |a|^{2-\beta} \Big\}.
\nonumber
\end{eqnarray}
This completes the proof.
\qed
\vskip 3mm

\noindent{\it Proof of Theorem \ref{Theorem:Infinite1}.} 
Note that $\xi\cap [-L,L]$, $L>0$ and $\xi$ 
satisfy ({\bf C.1}) and ({\bf C.2})
with the same constants $C_0, C_1, C_2$ and indices $\alpha, \beta$.
By virtue of Lemma \ref{lemma:4_2_5}
we see that there exists $C_3 >0$ such that
$$
|\Phi (\xi\cap [-L,L],a, iy)|
\le \exp \Big[ C_3 \Big\{|y|^{\theta} 
+ (|a|^{\theta} \vee 1) \Big\} \Big],
$$
$\forall L>0, \; \forall a\in \supp \xi, \forall y \in \R$.
Since for any $y\in\R$
$$
\Phi (\xi\cap [-L,L],a,iy) \to \Phi (\xi,a,iy), \quad L\to\infty,
$$
we can apply Lebesgue's convergence theorem 
to (\ref{eqn:KN1}) and obtain
$$
\lim_{L\to\infty}\mbK^{\xi\cap [-L,L]}\left(s, x; t, y \right)
=\mbK^{\xi}\left(s, x; t, y\right).
$$
Since for any $(s,t) \in (0, \infty)^{2}$ and 
any compact interval $I \subset \R$
$$
\sup_{x, y \in I} \Big|
\mbK^{\xi \cap [-L, L]}(s, x; t, y) \Big| < \infty,
$$
we can obtain the convergence of generating functions
for multitime correlation functions;
${\cal G}^{\xi \cap [-L, L]}[\chi] \to
{\cal G}^{\xi}[\chi]$ as $L \to \infty$.
It implies 
$\P_{\xi \cap [-L, L]} \to \P_{\xi}$ as $L \to \infty$
in the sense of finite dimensional distributions
and the proof is completed.
\qed

\subsection{Proofs of (\ref{eqn:Jacobi_kernel})
and (\ref{eqn:conv})}

\noindent{\it Proof of (\ref{eqn:Jacobi_kernel}).} 
Since $\xi^{\Z}=\eta^{1} \in \mX_0$,
we can start from the expression
of the correlation kernel (\ref{eqn:K1}) in 
Theorem \ref{Theorem:Infinite1}.
For $\ell \in \Z$, $z \in \C$
\begin{eqnarray}
\Phi(\xi^{\Z}, \ell, z) &=&
\prod_{j \in \Z, j \not=\ell}
\left(1- \frac{z-\ell}{j-\ell} \right) \nonumber\\
&=& \frac{\sin \{\pi(z-\ell)\}}{\pi(z-\ell)} \nonumber\\
&=& \frac{1}{2\pi} \int_{|k| \leq \pi} dk \,
e^{i k (z-\ell)},
\nonumber
\end{eqnarray}
since
$
\prod_{n\in\N}
(1-x^2/n^2)
= \sin (\pi x)/(\pi x)
$.
Then
\begin{equation}
\mbK^{\xi^{\Z}}(s, x; t, y)+{\bf 1}(s>t) p(s-t, x|y)
= \sum_{\ell\in\Z} p(s, x |\ell) I(t, y, \ell),
\label{eqn:eqA}
\end{equation}
where
\begin{eqnarray}
I(t, y,\ell) &=& \int_{\R} d y' \,
p(t, -iy | y')
\frac{1}{2\pi} \int_{|k| \leq \pi} dk \,
e^{i k (iy'-\ell)}
\nonumber\\
&=& \frac{1}{2 \pi} \int_{|k| \leq \pi} dk \,
e^{k^2 t/2 + i k (y-\ell)}.
\nonumber
\end{eqnarray}
By definition (\ref{eqn:theta3}) of $\vartheta_3$,
we can rewrite (\ref{eqn:eqA}) as
\begin{eqnarray}
&& \frac{1}{2 \pi} \int_{|k| \leq \pi} dk \,
e^{k^2 (t-s)/2+ i k (y-x)} \nonumber\\
&& \qquad \times
\vartheta_3 \left( \frac{1}{2 \pi i s} (x-i k s),
- \frac{1}{2 \pi i s} \right)
e^{-\pi i (x-i k s)^2/(2 \pi i s)}
\sqrt{\frac{i}{2 \pi i s}}.
\nonumber
\end{eqnarray}
Use the functional equation satisfied by
$\vartheta_3(v, \tau)$ 
(see, for example, Section 10.12 in \cite{AAR99}),
$$
\vartheta_3(v,\tau)
= \vartheta_3 \left( \frac{v}{\tau}, 
-\frac{1}{\tau} \right) e^{-\pi i v^2/\tau}
\sqrt{\frac{i}{\tau}},
$$
and the integral representation of the
heat kernel (\ref{eqn:p1}).
Then 
(\ref{eqn:Jacobi_kernel}) is obtained. \qed
\vskip 0.5cm

\noindent{\it Proof of (\ref{eqn:conv}).} 
By the definition (\ref{eqn:theta3}) of $\vartheta_3$,
for $s, t, u > 0$
\begin{eqnarray}
&& \mbK^{\xi^{\Z}}(u+s, x; u+t, y)
-{\bf K}_{\sin}(t-s, y-x) \nonumber\\
&=& \frac{e^{-2 \pi i x}}{2 \pi} 
\int_{|k| \leq \pi} dk \, 
e^{k^2(t-s)/2+ i k (y-x) 
- 2 \pi (u+s)(\pi+k)} \nonumber\\
&& + \frac{e^{2 \pi i x}}{2 \pi}
\int_{|k| \leq \pi} dk \, 
e^{k^2(t-s)/2+ i k (y-x) 
- 2 \pi (u+s)(\pi-k)} \nonumber\\
&& + \sum_{\ell \in \Z \setminus \{-1,0,1\}}
\frac{e^{2 \pi i x\ell}}{2 \pi}
\int_{|k| \leq \pi} dk \, 
e^{k^2 (t-s)/2+ i k (y-x) 
- 2 \pi (u+s) \ell(\ell\pi - k)}.
\nonumber
\end{eqnarray}
Then we see for any $u > 0$
\begin{eqnarray}
&& \Big| \mbK^{\xi^{\Z}}(u+s, x; u+t, y)
-{\bf K}_{\sin}(t-s, y-x) \Big|
\nonumber\\
&\leq& \Big( e^{\pi^2 (t-s)/2} \vee 1 \Big)
\left\{ \frac{1}{\pi} \int_{|k| \leq \pi} dk \, 
e^{-2 \pi (u+s)(\pi+k)}
+ 2 \sum_{\ell \geq 2} e^{-2 \pi^2 (u+s) \ell} \right\}
\nonumber\\
&=& \Big( e^{\pi^2 (t-s)/2} \vee 1 \Big)
\left\{ \frac{1-e^{-4 \pi^2 (u+s)}}{2\pi^2 (u+s)}
+ \frac{2 e^{-4 \pi^2 (u+s)}}{1-e^{-2 \pi^2 (u+s)}} 
\right\} \leq \frac{C}{u},
\nonumber
\end{eqnarray}
where $C>0$ depends on $t$ and $s$, but
does not on $u$.
This completes the proof of (\ref{eqn:conv}). \qed

\vskip 3mm
\noindent{\bf Remark} 
Since this relaxation process
$(\P_{\xi^{\Z}}, \Xi(t), t \in [0, \infty))$ is
determinantal with $\mbK^{\xi^{\Z}}$,
at any intermediate time $0 < t < \infty$
the particle distribution on $\R$ is
the determinantal point process with the
spatial correlation kernel
$\mbK^{\xi^{\Z}}(x,t; y, t), x, y \in \R$.
It should be noted that this spatial correlation kernel
is not symmetric,
$$
\mbK^{\xi^{\Z}}(t, x; t, y)
= \sum_{\ell \in \Z}
e^{2 \pi i x \ell-2 \pi^2 t \ell^2}
\frac{\sin \Big[ \pi 
\{(y-x)-2 \pi i t \ell\} \Big]}
{\pi \{(y-x)-2 \pi i t \ell\} },
$$
$x, y \in \R, 0 < t < \infty$.

\subsection{Proof of Theorem \ref{Theorem:Infinite2}}

In this subsection we prove Theorem \ref{Theorem:Infinite2}.
Suppose $\xi\in \mY^{\kappa}_m \subset \mY$. 
For $k\in\Z$ we can take $\lb_k$ and $\ub_k$ such that
$\lb_{-k-1}=-\ub_k, \ub_{-k-1}=-\lb_k$,
\begin{eqnarray}
&&\left[\lb_k, \ub_k \right] \subset (g^\kappa(k), g^\kappa(k+1)),
\quad
\ub_k-\lb_k \ge \frac{g^\kappa(k+1)-g^\kappa(k)}{2m(\xi,\kappa)+1}, 
\nonumber\\
&&\xi\big(\left[\lb_k, \ub_k \right]\big)=0
\quad
\mbox{and}
\quad
\xi\big(\big\{ (\ub_{k-1}+\lb_{k})/2 \big\}\big) = 0.
\nonumber
\end{eqnarray}
We put 
$I_k = \left[\lb_k, \ub_k \right]$,
$\varepsilon_k = |I_k| = \ub_k-\lb_k$,
$c_k = (\ub_{k-1}+\lb_{k})/2$, and
$\Delta_k = (\lb_{k}-\ub_{k-1})/2$.
Note that $[\ub_{-k-1}, \lb_{-k}]=-[\ub_{k-1}, \lb_{k}]$,
$I_{-k-1}=-I_{k}$, $\varepsilon_{-k-1}=\varepsilon_{k},
k \in \N_0$.
Then we define the $k$-th {\it cluster} 
in the configuration $\xi$ by
$$
\mC_k =\xi \cap [\ub_{k-1},\lb_{k}].
$$
It is easy to see that 
$\sum_{k \in\Z}\mC_{k}=\xi$,
and for each $k\in\Z$
\begin{eqnarray}
&& {|\mC_k|} \equiv \mC_k (\R) = \xi([\ub_{k-1},\lb_{k}]) 
\le 2m(\xi, \kappa), 
\label{Ck:1}
\\
&& |x-y| \ge \varepsilon_{k-1}\wedge\varepsilon_{k},
\quad x \in \supp \mC_k, \quad
y \in \supp (\xi-\mC_k).
\label{Ck:2}
\end{eqnarray}
Let $\v_k=(v_{k\ell})_{\ell=1}^{{|\mC_k|}}$
be the increasing sequence with 
$\sum_{\ell=1}^{{|\mC_k|}} \delta_{v_{k\ell}}=\mC_k$.
See Figure 1.
For $ a\in\supp\xi$, we denote by $\mC^a$ 
the cluster containing $a$.
Remark that, when $\xi_n$ converges to $\xi$ vaguely as $n\to\infty$,
we can take the $k$-th cluster
$\mC_k(\xi_n)$ of $\xi_n$ so that
it converges to the $k$-th cluster $\mC_k(\xi)$
of $\xi$ vaguely as $n\to\infty$.

\vskip 3mm
\begin{figure}[ch]
\begin{center}
{\footnotesize
\input Fig_Katori_Tanemura_08.tex
}
\end{center}
\caption{The clusters}
\end{figure}

We introduce $\C$-valued functions 
$\Psi_k(t,\xi,z, x)$, $k\in\Z$,
$t \geq 0$, $\xi\in\mY$, $z\in\C$, $x \in \R$,
\begin{eqnarray}
&&\Psi_k(t,\xi,z, x)=\Phi(\xi- \mC_k,{c_k},z)
\sum_{\ell=1}^{{|\mC_k|}}(z-{c_k})^{\ell-1}
(-1)^{{|\mC_k|}-\ell-1}
\nonumber\\
&&\qquad
\times \left\{ \Theta_{k,\ell}(t,\xi,x)+
\sum_{q={|\mC_k|}}^\infty \Theta_{k,q}(t,\xi,x)
s_{(q-{|\mC_k|}|{|\mC_k|}-\ell-1)}(\v_k-{c_k}) \right\},
\label{PSI2}
\end{eqnarray}
if $|\mC_k|\not=0$, and $\Psi_k(t,\xi,z)=0$, otherwise, where
$s_{(k|\ell)}$ is the Schur function associated with
the partition $(k|\ell)$ in Frobenius' notation, and
\begin{equation}
\Theta_{k,q}(t,\xi,x)=
\sum_{r=0}^q \frac{1}{(q-r)!}
\left( - \frac{1}{\sqrt{2t}} \right)^{q-r}
H_{q-r}\left(\frac{c_k-x}{\sqrt{2t}} \right)
h_r \left(
\left(\frac{1}{u-{c_k}} \right)_{u\in \xi-\mC_k}
\right)
\label{eqn:Theta}
\end{equation}
with the Hermite polynomials $H_k, k\in\N_0$,
and with the complete symmetric functions $h_k, k\in\N_0$.

\begin{lem}\label{lemma:4_3_2}
Suppose that $\xi\in \mY_0$.
Then for $k \in \Z, t \geq 0, x \in \R, z \in \C$
$$
\int_\R \mC_k(dx') e^{-(x'-x)^2/(2t)}\Phi(\xi,x',z)=
e^{-({c_k}-x)^2/(2t)}\Psi_k(t,\xi, z,x).
$$
\end{lem}
\vskip 3mm

\noindent {\it Proof.} 
From definitions of $\mC_k, k\in\Z$ and $\Phi$,
we have
\begin{eqnarray}
&&\int_\R \mC_k(dx') e^{-(x'-x)^2/(2t)}\Phi(\xi,x',z)
\nonumber\\
&&=\int_\R \mC_k(dx') e^{-(x'-x)^2/(2t)} 
\prod_{u\in \xi - \mC_k} \frac{z -u}{x'-u}
\prod_{v\in \mC_k - \delta_{x'}} \frac{z -v}{x'-v}
\nonumber
\\
&& = e^{-{(c_k-x)}^2/(2t)}\int_\R \mC_k(dx')
e^{-(x'-{c_k})(x'+{c_k}-2x)/(2t)}
\nonumber\\
&&\qquad\times
\prod_{u\in \xi - \mC_k}
\frac{(z-{c_k}) -(u-{c_k})}{(x'-{c_k})-(u-{c_k})}
\prod_{v\in \mC_k-\delta_{x'}}
\frac{(z-{c_k}) -(v-{c_k})}{(x'-{c_k})-(v-{c_k})}
\nonumber\\
&& =
e^{-{(c_k-x)}^2/(2t)}
\sum_{j=1}^{{|\mC_k|}}
\psi_k(t,\xi, v_{kj}-{c_k},z,x)
\frac{a_\delta(\v_k -{c_k};j;z-{c_k})}{a_\delta(\v_k-{c_k})},
\nonumber
\end{eqnarray}
where
$$
a_\delta(\x_m;j;y)=a_{\delta}(x_1,\dots,x_{j-1},y,x_{j+1},\dots,x_m),
$$
and
$$
\psi_k(t,\xi,x',z,x) = 
\Phi(\xi- \mC_k, x'+{c_k},z)
\exp\left( -\frac{2(c_k-x) x'+{x'}^2}{2t} \right).
$$
Now we introduce $\widetilde{\Theta}_{k,q}$'s
as the coefficients of the expansion
$$
\psi_k(t,\xi,x',z,x)= \sum_{q \in \N_0}
\widetilde{\Theta}_{k,q}(t,\xi,z,x) {x'}^q.
$$
Then we have
\begin{eqnarray}
&&\frac{1}{a_\delta(\v_k-{c_k})}\sum_{j=1}^{{|\mC_k|}}
\psi_k(t,\xi,v_{kj}-{c_k},z,x)a_\delta(\v_k -{c_k};j;z-{c_k})
\nonumber\\
&& \, =\frac{1}{a_\delta (\v_k-{c_k})}
\sum_{\ell=1}^{{|\mC_k|}}(z-{c_k})^{\ell-1}(-1)^{{|\mC_k|}-\ell-1}
\nonumber
\\
&& \, \times
\det
\left( 
\begin{array}{cccc}
\psi_k(t,\xi,v_{k1}-{c_k},z,x)&\psi_k(t,\xi,v_{k2}-{c_k},z,x)
&\cdots&\psi_k(t,\xi,v_{k{|\mC_k|}}-{c_k},z,x)
\\
(v_{k1}-{c_k})^{{|\mC_k|}-1}&(v_{k2}-{c_k})^{{|\mC_k|}-1}&\cdots&
(v_{k{|\mC_k|}}-{c_k})^{{|\mC_k|}-1}
\\
\dots&\dots&\dots&\dots
\\
(v_{k1}-{c_k})^{\ell+1}&(v_{k2}-{c_k})^{\ell+1}&\cdots&
(v_{k{|\mC_k|}}-{c_k})^{\ell+1}
\\
(v_{k1}-{c_k})^{\ell-1}&(v_{k2}-{c_k})^{\ell-1}&\cdots&
(v_{k{|\mC_k|}}-{c_k})^{\ell-1}
\\
\dots&\dots&\dots&\dots
\\
v_{k1}-{c_k}&v_{k2}-{c_k}&\cdots&v_{k {|\mC_k|}}-{c_k}
\\
1&1&\dots&1
\end{array}
\right)
\nonumber\\
&&\, = \sum_{\ell=1}^{{|\mC_k|}}(z-{c_k})^{\ell-1}
(-1)^{{|\mC_k|}-\ell-1}
\nonumber\\
&&\qquad\qquad \times \left\{ \widetilde{\Theta}_{k,\ell}(t,\xi,z,x)+
\sum_{q={|\mC_k|}}^\infty \widetilde{\Theta}_{k,q}(t,\xi,z,x)
s_{(q-{|\mC_k|}|{|\mC_k|}-\ell-1)}(\v_k-{c_k}) \right\}.
\nonumber
\end{eqnarray}
Then, to prove the lemma, 
it is enough to show the equality
\begin{equation}
\widetilde{\Theta}_{k,q}(t,\xi,z,x)
= \Phi(\xi- \mC_k,{c_k},z)
\Theta_{k,q}(t,\xi,x),
\quad t \geq 0, \; \xi\in\mY_0, \; z\in\C, \; x \in \R,
\label{Psik}
\end{equation}
for $|\mC_k| \not= 0$.
From the formula (\ref{complete}), we have
\begin{eqnarray}
&&\Phi(\xi - \mC_k, x'+{c_k},z) 
= \prod_{u \in \xi-\mC_k}\frac{z-u}{x'-(u-{c_k})}
\nonumber\\
&&\qquad
= \prod_{u \in \xi- \mC_k}\frac{u-z}{u-{c_k}}
\prod_{u \in \xi-  \mC_k}\frac{1}{1-x'/(u-{c_k})}
\nonumber\\
&&\qquad=\prod_{u \in \xi- \mC_k}\frac{u-z}{u-{c_k}}
\sum_{r \in \N_0}h_r \left(
\left(\frac{1}{u-{c_k}} \right)_{u \in \xi- \mC_k}
\right){x'}^r
\nonumber\\
&&\qquad= \Phi(\xi- \mC_k,{c_k},z)
\sum_{r \in \N_0}h_r \left(
\left(\frac{1}{u-{c_k}} \right)_{u \in \xi- \mC_k}
\right){x'}^r.
\label{CSF}
\end{eqnarray}
By the formula (\ref{eqn:Hermite3}), 
we have
\begin{equation}
\exp \left( - \frac{2(c_k-x) x'+{x'}^2}{2t} \right)
= \sum_{k \in \N_0} \frac{1}{k!}
\left(- \frac{x'}{\sqrt{2t}} \right)^{k}
H_k \left( \frac{c_k-x}{\sqrt{2t}} \right).
\label{GFH}
\end{equation}
Combining (\ref{CSF}) and (\ref{GFH}),
we have
\begin{eqnarray}
&&\psi_k(t,\xi,x',z,x)
\nonumber\\
&&=\Phi(\xi- \mC_k,{c_k},z)
\sum_{r \in \N_0}
h_r \left(
\left(\frac{1}{u-{c_k}} \right)_{u \in \xi- \mC_k}
\right){x'}^r \nonumber\\
&& \hskip 4cm \times
\sum_{k \in \N_0} \frac{1}{k!}
\left( - \frac{x'}{\sqrt{2t}} \right)^{k}
H_k \left( \frac{c_k-x}{\sqrt{2t}} \right)
\nonumber\\
&&=\Phi(\xi- \mC_k,{c_k},z)
\sum_{q \in \N_0} {x'}^q
\sum_{r=0}^q \frac{1}{(q-r)!}
\left( -\frac{1}{\sqrt{2t}} \right)^{q-r}
\nonumber\\
&& \hskip 4cm \times
H_{q-r} \left(\frac{c_k-x}{\sqrt{2t}} \right)
h_r \left(
\left(\frac{1}{u-{c_k}} \right)_{u \in \xi-\mC_k}
\right).
\nonumber
\end{eqnarray}
Then, by definition (\ref{eqn:Theta}),
(\ref{Psik}) is proved.
\qed
\vskip 3mm
\begin{lem}\label{lemma:4_3_1}
Assume that {\rm ({\bf C.3})} holds 
with some $\kappa\in (1/2,1)$ and $m\in\N$. 

\noindent{\rm (i)} 
Suppose that $\alpha\in (1/\kappa,2)$.
Then there exists $C_4(\kappa,m,\alpha) >0$ such that
\begin{equation}
M_\alpha\left(\tau_{-a}(\xi- \mC^a)\right)
\le C_4(\kappa,m) (|a|\vee 1)^{(1-\kappa)/\kappa}
\quad \forall a \in \supp \xi,
\label{estimate:C_3_1}
\end{equation}
and {\rm ({\bf C.2}) (i)} holds, that is, 
there exists $C_1=C_1(\alpha,\xi)$ such that
\begin{equation}
M_{\alpha}(\xi)\le C_1.
\label{estimate:C_3_2}
\end{equation}

\noindent{\rm (ii)} 
Suppose that $\beta \in (0,2\kappa-1)$.
Then $\xi-\mC^a-\widehat{\mC^a}$ satisfies
{\rm ({\bf C.2}) (ii)} $\forall a \in \supp \xi$,
where $\widehat{\mC^a}=\mC_{-k}$ in case ${\mC^a}=\mC_{k}$.
That is, 
there exists $C_2(\kappa,m) >0$ such that
\begin{equation}
M_1\left( \tau_{-a^2}(\xi- \mC^a-\widehat{\mC^a})^{\langle 2 \rangle}
\right)
\le C_2(\kappa,m) (|a|\vee 1)^{-\beta}
\quad \forall a \in \supp \xi.
\label{estimate:C_3_3}
\end{equation}
\end{lem}
\vskip 3mm

\noindent {\it Proof.} 
By simple calculations 
we see that there exists a positive constant
$C(\kappa)$ such that
\begin{equation}
\label{bkappa}
M_\alpha (\tau_{-a}\eta^{\kappa})
\le C(\kappa) (|a|\vee 1)^{(1-\kappa)/\kappa}
\quad \forall a\in \supp\eta^{\kappa}.
\end{equation}
Suppose that $\mC^a=\mC_k$, $k\in \Z$.
Then $\xi - \mC^a= \xi \cap [\ub_{k-1},\lb_{k}]^{\rm c}$. 
We divide the set $[\ub_{k-1},\lb_{k}]^{\rm c}$
into the following four sets:
\begin{eqnarray}
&& A_1 = \Big( -\infty, g^{\kappa}(k-2) \Big], 
\qquad
A_2 = \Big( g^{\kappa}(k-2), \ub_{k-1} \Big),
\nonumber\\
&&
A_3 = \Big( \lb_{k}, g^{\kappa}(k+2) \Big),
\qquad
A_4 = \Big[ g^{\kappa}(k+2), -\infty \Big).
\nonumber
\end{eqnarray}
Then we have
$$
\left( 
\int_{\R} \frac{(\xi- \mC^a)(dx)}{|x-a|^\alpha}\right)^{1/\alpha}
\le \sum_{j=1}^4 \left(
\int_{A_j}\frac{\xi(dx)}{|x-a|^\alpha}
\right)^{1/\alpha}.
$$
From (\ref{Ck:1}) and (\ref{Ck:2}),
we have 
\begin{eqnarray}
\int_{A_1} \frac{\xi(dx)}{|x-a|^\alpha}
&\le& m \sum_{-\infty <\ell\le k-2}
\frac{1}{|g^{\kappa}(\ell)-g^{\kappa}(k-1)|^\alpha},
\nonumber\\
\int_{A_2} \frac{\xi(dx)}{|x-a|^\alpha}
&\le& 2m
\left(\frac{1}{\varepsilon_{k-1}}\right)^\alpha,
\nonumber\\
\int_{A_3} \frac{\xi(dx)}{|x-a|^\alpha}
&\le& 2m
\left(\frac{1}{\varepsilon_{k}}\right)^\alpha,
\nonumber\\
\int_{A_4} \frac{\xi(dx)}{|x-a|^\alpha}
&\le& m \sum_{k+2 \le \ell< \infty}
\frac{1}{|g^{\kappa}(\ell)-g^{\kappa}(k+1)|^\alpha}.
\nonumber
\end{eqnarray}
Combining these estimates with (\ref{bkappa}), we have
$$
\left( 
\int_{\R} \frac{(\xi- \mC^a)(dx)}{|x-a|^\alpha}
\right)^{1/\alpha}
\le {\cal O}
\Big((|g^{\kappa}(k-1)|\vee 
|g^{\kappa}(k+1)|\vee 1)^{(1-\kappa)/\kappa} \Big),
\quad |k| \to\infty.
$$
Since $\max_{k-1 \le j \le k+1} |g^{\kappa}(j)|
\le 2(|a|\vee 1)$,
we obtain (\ref{estimate:C_3_1}).
The estimate (\ref{estimate:C_3_2}) is derived from
(\ref{estimate:C_3_1}) with $a=0$ and $\mC^a=\mC_0$,
and the fact that $M_{\alpha}(\mC_0)<\infty$.
Noting that $(\xi- \mC^a-\widehat{\mC^a})^{\langle 2 \rangle}$
satisfies {\rm ({\bf C.3})} 
with $2\kappa$ and $2m$,
we obtain (\ref{estimate:C_3_3})
by a similar argument given above 
to show (\ref{estimate:C_3_1}).
This completes the proof.
\qed
\vskip 3mm

\begin{lem}\label{lemma:4_2_4}
Let $\alpha\in (1,2)$ and $|a| \ge 1$.
Assume that {\rm ({\bf C.1})} and the condition
\begin{equation}
M_\alpha(\tau_{-a}\xi) \le C_5 |a|^\gamma
\label{old_C2}
\end{equation}
with some $\gamma>0$ and $C_5>0$ are satisfied. 
Then there exists $C_6=C_6(\alpha, \beta, C_1, C_5) >0$ 
such that
$$
|M({\tau_{-a}\xi})-M(\xi)|\le C_6 |a|^{\delta_1},
$$
where $\delta_1=\alpha(1+\gamma)-1$.
\end{lem}

\noindent {\it Proof.}
From Lemma \ref{lemma:4_2_2} and the fact that
$M_1({\tau_{-a}\xi},L)$ is increasing in $L$, we see that
$$
\max_{0\le L \le L_0}M_1({\tau_{-a}\xi},L)= M_1({\tau_{-a}\xi},L_0)
\le (2M_\alpha({\tau_{-a}\xi}))^{\alpha{\delta_1}/({\delta_1}-\alpha+1)}
\le C |a|^{\delta_1}
$$
from (\ref{old_C2}) with a constant $C>0$.
Combining this estimate with Lemma \ref{lemma:4_2_2},
we have
\begin{equation}
M_1({\tau_{-a}\xi},L)
\le C |a|^{\delta_1} \vee L^{\delta_1}.
\label{old_lemma423}
\end{equation}
We assume $a\not= 0$.
By the definitions of $M(\xi)$ and $M({\tau_{-a}\xi})$,
$$
|M({\tau_{-a}\xi})-M(\xi)|
\le \frac{1+\xi(\{0\})}{|a|}
+ |a|\int_{\{a,0\}^{\rm c}} \frac{\xi(dx)}{|x(x-a)|}.
$$
We divide the set $\{a,0 \}^{\rm c}$ into 
the three disjoint subsets 
$\{x:0<|x|< 2|a|, 2|a-x|>|a|\}$,
$\{x:|x|\ge 2|a|\}$
and $\{x:0<|x|< 2|a|, 0<2|a-x|\le |a| \}$.
By simple calculation, we see
$$
\int_{0<|x|< 2|a|, 2|a-x|>|a|}
\frac{\xi(dx)}{|x(x-a)|}
\le \frac{2}{|a|}
\int_{0<|x|<2|a|}\frac{\xi(dx)}{|x|}
=\frac{2}{|a|}M_1(\xi,2|a|).
$$
Since $|x-a|\ge |x|-|a|\ge |x|/2$,
if $|x|\ge 2|a|$,
$$
\int_{|x|\ge 2|a|}\frac{\xi(dx)}{|x(x-a)|}
\le 2\int_{|x|\ge 2|a|}\frac{\xi(dx)}{|x|^2}
\le 2^{\alpha-1}M_\alpha(\xi)^\alpha |a|^{\alpha-2}.
$$
Since $|x|\ge |a|-|a-x|\ge |a|/2$,
if $2|a-x|\le |a|$,
$$
\int_{0<|x|< 2|a|, 0<2|a-x|\le |a|}\frac{\xi(dx)}{|x(x-a)|}
\le\frac{2}{|a|}
\int_{0<2|a-x|\le|a|}\frac{\xi(dx)}{|x-a|}
=\frac{2}{|a|}M_1 \left({\tau_{-a}\xi},\frac{|a|}{2} \right).
$$
Combining the above estimates with the fact $|a|^{-1} \le 1$,
we have 
$$
|M({\tau_{-a}\xi})-M(\xi)|\le 2^{\alpha-1}M_\alpha(\xi)^\alpha |a|^{\alpha-1}
+2M_1(\xi, 2|a|)+2M_1
\left({\tau_{-a}\xi},\frac{|a|}{2} \right)+2.
$$
Then the lemma is derived from 
(\ref{old_C2}) and (\ref{old_lemma423}). \qed

\vskip 3mm

The following is a key lemma to prove
Theorem \ref{Theorem:Infinite2}.

\begin{lem}
\label{lemma:4_3_3}
Let $t \geq 0, x \in \R$,
$\xi\in \mY^{\kappa}_{m} \subset \mY$ 
with $\kappa\in (1/2,1)$ and $m\in\N$.
Then for any $\theta\in (3-2\kappa, 2)$ there exist positive constants 
$C_7= C_7(t, \kappa, C_0,x)$ and
$\widehat{C}_7= \widehat{C}_7(t, \kappa, m, \theta, C_0,x)$
such that
$$
|\Psi_k (t, \xi, iy,x)|
\le \widehat{C}_7 \exp \Big[ C_7 \Big\{|y|^{\theta} 
+ |{c_k}|^{\theta} \Big\} \Big],
\quad \forall y\in\R, \; \forall k\in \Z.
$$
\end{lem}

\noindent {\it Proof.}
We note the equality
$$
\Phi(\xi-\mC_k, c_k,iy)=
\Phi(\xi-\mC_k-\mC_{-k}, c_k,iy)
\Phi(\mC_{-k}, c_k,iy).
$$
Let $\beta\in (0,2\kappa-1)$ and $\alpha=(1/\kappa,2)$.
By virtue of Lemma \ref{lemma:4_3_1},
we can apply Lemma \ref{lemma:4_2_5}
for $\xi- \mC_k-\mC_{-k}$
and see that there exist positive constant $C_3$ 
and $\theta \in (3-2\kappa,2)$ such that
$$
|\Phi (\xi-\mC_k-\mC_{-k}, c_k,iy)|
\le \exp \Big[ C_3 \Big\{|y|^{{\theta}} 
+ |{c_k}|^{{\theta}} \Big\} \Big],
\quad y\in\R, \; k\in \Z.
$$
Here we used the fact that 
$3-2\kappa>1/\kappa$ for $\kappa\in (1/2,1)$.
Since $\Phi(\mC_{-k}, c_k,iy)$ is 
a polynomial function of $y$, we have
$$
|\Phi (\xi-\mC_k, c_k,iy)|
\le \widehat{C}_3\exp \Big[ C_3 \Big\{|y|^{{\theta}} 
+ |{c_k}|^{{\theta}} \Big\} \Big],
\quad y\in\R, \; k\in \Z,
$$
for some $\widehat{C}_3>0$.
Hence, from the definition (\ref{PSI2}) of $\Psi_k (t,\xi, z, x)$,
to prove the lemma
it is enough to show the following estimates:
for any $\ell=1,2,\dots, {|\mC_k|}$,
\begin{eqnarray}
&&|(z-{c_k})^{\ell-1}| 
= {\cal O} (|z|^{{|\mC_k|}}\vee |{c_k}|^{{|\mC_k|}}),
\quad |k| \to\infty, \; |z|\to\infty,
\label{est_1}
\\
&&|\Theta_{k,\ell}(t,\xi,x)|
= {\cal O} ( |{c_k}|^\ell),
\quad |k| \to\infty,
\label{est_2}
\\
\label{est_0}
&&\sum_{q={|\mC_k|}}^\infty \Theta_{k,q}(t,\xi,x)
s_{(q-{|\mC_k|}|{|\mC_k|}-\ell-1)}(\v_k-{c_k})
\le \exp\Big[ C (|{c_k}|^{\theta'} \vee 1) \Big],
\quad k\in \Z, \qquad
\end{eqnarray}
with some $C=C(t,x) >0$ and $\theta'<\theta$. Since
(\ref{est_1}) and (\ref{est_2}) can be confirmed easily,
here we show only the proof of (\ref{est_0}).
Since
$|v_{k,\ell} -{c_k}|\le \Delta_k$,
$1\le \ell \le {|\mC_k|}$,
from the fact (\ref{skl1})
$$
s_{(q-{|\mC_k|}|{|\mC_k|}-\ell-1)}(\v_k-{c_k})
\le {q-\ell-1 \choose {|\mC_k|}-\ell-1} {q \choose \ell} \Delta_k^q
\le q^{{|\mC_k|}} \Delta_k^q,
\quad q\in\N.
$$
Put
$\overline{\Delta}_k =
\Delta_k+ (\varepsilon_{k-1}\wedge\varepsilon_{k})/2$,
and remind that 
$\overline{\Delta}_k= {\cal O}(c_k^{(\kappa-1)/\kappa})$,
$|k| \to\infty$.
Then we have
$$
s_{(q-{|\mC_k|}|{|\mC_k|}-\ell-1)}(\v_k-{c_k})
\le C' {{\overline{\Delta}_k}}^q,
\quad k\in\Z, \; q\in\N, 
$$
with some positive constant $C' >0$.
Then
\begin{eqnarray}
&&\Theta_{k,q}(t,\xi,x)
s_{(q-{|\mC_k|}|{|\mC_k|}-\ell-1)}(\v_k-{c_k})
\nonumber\\
&&\quad\le C' \sum_{r=0}^q 
\frac{1}{(q-r)!}
\left(\frac{\overline{\Delta}_k}{\sqrt{2t}}\right)^{q-r}
\left|H_{q-r}\left( \frac{c_k-x}{\sqrt{2t}} \right) \right|
{\overline{\Delta}_k}^r
\left|h_r \left(
\left(\frac{1}{u-{c_k}} \right)_{u \in \xi- \mC_k}
\right)\right|,
\nonumber
\end{eqnarray}
and thus
\begin{eqnarray}
&&\sum_{q={|\mC_k|}}^\infty \Theta_{k,q}(t,\xi,x)
s_{(q-{|\mC_k|}|{|\mC_k|}-\ell-1)}(\v_k-{c_k})
\nonumber\\
&&\quad
\le  C' \sum_{q \in \N_0}
\frac{1}{q!}
\left(\frac{\overline{\Delta}_k}{\sqrt{2t}}\right)^{q}
\left|H_{q}\left( \frac{c_k-x}{\sqrt{2t}} \right) \right|
\sum_{r \in \N_0}
{\overline{\Delta}_k}^r
\left|h_r \left(
\left(\frac{1}{u-{c_k}} \right)_{u \in \xi- \mC_k}
\right)\right|.
\nonumber
\end{eqnarray}
Since
$$
\left| \left.
\frac{d^k}{dz^k} e^{2zy-z^2}
\right|_{z=0} \right|
\le \left.
\frac{d^k}{dz^k} e^{2z|y|+z^2} \right|_{z=0},
\quad k\in\N, y \in \R, 
$$
we obtain from (\ref{GFH})
\begin{eqnarray}
&&\sum_{q\in\N_0}
\frac{1}{q!}
\left(\frac{\overline{\Delta}_k}{\sqrt{2t}}\right)^{q}
\left|H_{q}\left( \frac{c_k-x}{\sqrt{2t}} \right) \right|
\nonumber\\
\label{eqt_1}
&& \quad 
\le \exp \left( \frac{2\overline{\Delta}_k(c_k-x)+
\overline{\Delta}_k^2}{2t} \right)
={\cal O}\Bigg(\exp \Big( 
\widetilde{C} {c_k}^{1+(\kappa-1)/\kappa}
\Big) \Bigg),
\end{eqnarray}
$|k| \to \infty$, 
with a constant $\widetilde{C}=\widetilde{C}(t,x)$. 
And if $(\xi- \mC_k)(u)\ge 1$, then
$|u-{c_k}|\ge \Delta_k+\varepsilon_{k-1} \wedge \varepsilon_{k}$
and
$$
\frac{1}{1-\overline{\Delta}_k /|u-{c_k}|}\le C m
$$
with a positive constant $C$.
Hence from (\ref{complete2}) 
\begin{eqnarray}
&&\sum_{r\in\N_0}
{\overline{\Delta}_k}^r
\left|h_r \left(
\left(\frac{1}{u-{c_k}} \right)_{u \in \xi- \mC_k}
\right) \right|
\nonumber\\
&&\quad \le \exp \Bigg\{
\Big|M \Big(\tau_{-{c_k}}(\xi-\mC_k) \Big) \Big|\overline{\Delta}_k
+ Cm{\overline{\Delta}_k}^2 
M_2 \Big(\tau_{-{c_k}}(\xi -\mC_k) \Big)^2
\Bigg\}.
\label{eqt_2}
\end{eqnarray}
Using Lemmas \ref{lemma:4_3_1} and \ref{lemma:4_2_4},
we see that
$$
\bigg|M\big(\tau_{-{c_k}}(\xi- \mC_k)\big)\bigg|\overline{\Delta}_k
= {\cal O} \Big(|{c_k}|^{\delta_1 +(\kappa-1)/\kappa} \Big),
\quad |k| \to\infty,
$$
with any
$\delta_1>\{1+(1-\kappa)/\kappa\}/\kappa-1
=1/\kappa^2-1$,
and
$$
{\overline{\Delta}_k}^2 M_2\big(\tau_{-{c_k}}(\xi-\mC_k)\big)^2
= {\cal O} \Big(|{c_k}|^{\alpha(1-\kappa)/\kappa}
\overline{\Delta}_k^\alpha \Big)
= {\cal O}(1), \quad |k| \to\infty.
$$
Since $1/\kappa^2-1+(\kappa-1)/\kappa +1+(\kappa-1)/\kappa
=1/\kappa^2+2(\kappa-1)/\kappa < 3-2\kappa$, 
for $\kappa\in(1/2,1)$,
(\ref{est_0}) is derived from (\ref{eqt_1}) and  (\ref{eqt_2}).
This completes the proof.
\qed

\vskip 3mm

\noindent{\it Proof of Theorem \ref{Theorem:Infinite2}.} \,
(i) By Lemmas \ref{lemma:4_3_2} and \ref{lemma:4_3_3},
if $\xi \in \mY_0$,
\begin{eqnarray}
&& \int_{\R} \xi(dx') p(s,x|x')
\int_{\R} dy' \, p(t,-iy|y') \Phi(\xi, x', iy)
\nonumber\\
&=& \sum_{k \in \Z} p(s,x|c_k)
\int_{\R}dy' \, p(t, -iy|y') 
\Psi_{k}(t, \xi, iy', x).
\label{eqn:eqA1}
\end{eqnarray}
Since this equality holds even if we replace
$\xi$ by $\xi \cap [-L, L]$ for any $L >0$,
\begin{eqnarray}
&& \lim_{L \to \infty} \mbK^{\xi \cap [-L,L]}
(s, x; t, y)
+{\bf 1}(s>t) p(s-t,x|y)
\nonumber\\
&=& \lim_{L \to \infty} \int_{\R}
\xi \cap [-L, L](dx') p(s,x|x')
\int_{\R} dy' \, p(t, -iy|y') 
\Phi(\xi \cap [-L, L], x', iy)
\nonumber\\
&=& \lim_{L \to \infty} \sum_{k \in \Z}
p(s, x|c_k) \int_{\R}dy' \, p(t, -iy|y')
\Psi_k(t, \xi \cap [-L,L], iy', x).
\nonumber
\end{eqnarray}
By Lemma \ref{lemma:4_3_3}, we can apply Lebesgue's 
convergence theorem to show that the limit is
$$
\sum_{k \in \Z} p(s,x|c_k) \int_{\R} dy' \,
p(t, -iy|y') \Psi_k(t, \xi, iy', x).
$$
We can repeat the argument
in the proof of Theorem \ref{Theorem:Infinite1}
given at the end of Section 4.2.
Then if $\xi \in \mY_0$,
$(\P_{\xi}, \Xi(t), t \in [0, \infty))$ is well-defined
with the correlation kernel
\begin{eqnarray}
\mbK^{\xi}(s, x; t, y)
&=& \sum_{k\in \Z} p(s, x|c_k)
\int_{\R} d y' \, p(t, -iy|y')
\Psi_k(t, \xi, i y',x) \nonumber\\
&& -{\bf 1}(s>t)p(s-t,x|y).
\label{eqn:kxi}
\end{eqnarray}
It is equal to (\ref{eqn:K1}) of Theorem
\ref{Theorem:Infinite1} by the equality (\ref{eqn:eqA1}).
When $\xi \in \mY \setminus \mY_0$, (\ref{eqn:eqA1})
is not valid.
For any $L > 0$, however, the equality
\begin{eqnarray}
&& \mbK^{\xi \cap [-L,L]}
(s, x; t, y)
+{\bf 1}(s>t) p(s-t,x|y)
\nonumber\\
&=& \sum_{k \in \Z}
p(s, x|c_k) \int_{\R}dy' \, p(t, -iy|y')
\Psi_k(t, \xi \cap [-L,L], iy', x)
\nonumber
\end{eqnarray}
holds by the continuity with respect to the initial
configuration for Dyson's model with finite particles.
Then, again by Lemma \ref{lemma:4_3_3} with 
Lebesgue's convergence theorem, we will 
obtain the result (\ref{eqn:kxi}). \\
(ii) By the fact (\ref{con:moderate}) and the definition of $\Psi_k$, 
we see that for any $k\in\N$, $t \geq 0$, 
and $x,y' \in \R$ 
$$ 
\lim_{n\to\infty}\Psi_k(t,\xi_n,iy',x)= \Psi_k(t,\xi,iy',x).
$$
By using Lemma \ref{lemma:4_3_3} we see that, 
for fixed $t \geq 0, x \in \R$,
there exist ${\theta}\in (1,2)$, $C_7=C_7(t,x) >0$, and
$\widehat{C}_7=\widehat{C}_7(t,x) >0$ such that 
$$
|\Psi_k (t,\xi_n,iy',x)|
\le \widehat{C}_7 \exp \Big[ C_7 \Big\{ |y'|^{{\theta}}  
+ |{c_k}|^{{\theta}} \Big\} \Big], 
\quad k \in \Z, y' \in \R, n \in \N.
$$
Therefore, by applying Lebesgue's convergence theorem, 
we obtain the theorem.
\qed
\vskip 1cm

\begin{small}
\noindent{\it Acknowledgments.} 
The present authors would like to thank
T. Shirai and H. Spohn for useful comments
on the manuscript.
A part of the present work was done
during the participation of M.K.
in the ESI program ``Combinatorics and Statistical Physics"
(March and May in 2008).
M.K. expresses his gratitude for 
hospitality of the Erwin Schr\"odinger Institute 
(ESI) in Vienna
and for well-organization of the program
by M. Drmota and C. Krattenthaler.
M.K. is supported in part by
the Grant-in-Aid for Scientific Research (C)
(No.21540397) of Japan Society for
the Promotion of Science.
H.T. is supported in part by
the Grant-in-Aid for Scientific Research 
(KIBAN-C, No.19540114) of Japan Society for
the Promotion of Science.
\clearpage

\end{small}
\end{document}